
\magnification=\magstep1
\input amstex
\documentstyle{amsppt}
\def \R {\Bbb R}
\def \C {\Bbb C}
\def \N {\Bbb N}

\def \N {\Bbb N}

\def \Ri {\Cal R}

\def \D {\Cal D}

\def \h {\frak h}

\def \al {\alpha}
\def \la {\lambda}
\def \ph {\varphi}
\def \del {\delta}
\def \eps {\varepsilon}
\def \om {\omega}
\def \lan {\langle}
\def \ran {\rangle}
\def \proof {\demo {Proof}}
\def \endproof {\qed\enddemo}
\def \de {\partial} \def \p {\de}

\def \half{\frac12}
\def \thalf{\textstyle{\half}}
\def \inv{^{-1}}
\def \rinv#1{^{(#1)}}
\def \supp {\text{\rm supp\,}}
\def \dom {\text{\rm dom\,}}

\def \sgn {\text{\rm sgn\,}}

\topmatter
\title $L^p$-spectral multipliers for the Hodge Laplacian acting on
1-forms on the Heisenberg group 
\endtitle
\rightheadtext{Spectral multipliers for the Hodge Laplacian}
\leftheadtext{D. M\"uller, M. M. Peloso and F. Ricci}
\author Detlef M\"uller, Marco M. Peloso
  and Fulvio Ricci
\endauthor
\address \!\!\!\!\!\!\!\!Christian-Albrechts-Universit\"at zu Kiel,  
Mathematisches
Seminar, Ludewig- Meyn-Str. 4, D-24098 Kiel, Germany
\endaddress
\email mueller\@math.uni-kiel.de\endemail \urladdr
http://analysis.math.uni-kiel.de/mueller
\endurladdr
\address \!\!\!\!\!\!\!\!Dipartimento di Matematica, Politecnico di
Torino, Corso Duca degli Abruzzi 24, 10129 Torino, Italy \endaddress
\email peloso\@calvino.polito.it\endemail \urladdr
  http://calvino.polito.it/\~peloso/
\endurladdr
\address \!\!\!\!\!\!\!\!Scuola Normale Superiore, Piazza dei Cavalieri
7, 56126 Pisa, Italy \endaddress
\email fricci\@sns.it\endemail \urladdr
  http://www.math.sns.it/HomePages/Ricci/
\endurladdr
\abstract We prove that, if $\Delta_1$ is the Hodge Laplacian acting
on differential 1-forms on the $(2n+1)$-dimensional Heisenberg group,  
and if $m$ is a
Mihlin-H\"ormander multiplier on the positive half-line, with  
$L^2$-order of smoothness
greater than $n+\half$, then $m(\Delta_1)$ is $L^p$-bounded for  
$1<p<\infty$. Our
approach leads to an explicit description of the spectral decomposition  
of $\Delta_1$ on
the space of $L^2$-forms in terms of the spectral analysis of the  
sub-Laplacian $L$
and the central derivative $T$, acting on scalar-valued functions.  
\endabstract 
\thanks
This work has been supported by the the IHP network HARP ``Harmonic  
Analysis and Related problems'' of the European Union.
\endthanks 
\subjclass 43A80, 42B15
\endsubjclass 
\endtopmatter

\document

\head Introduction \endhead

The $(2n+1)$-dimensional Heisenberg group $H_n$ has a (unique modulo
dilations) left-invariant Riemannian structure which is invariant under  
the
action of the unitary group $U(n)$ by automorphisms (i.e. the natural  
action on the
$\C^n$-component, when $H_n$ is realized as $\C^n\times\R$). Various
differential-geometric aspects of this structure have been analyzed in
the literature \cite{DT, L, R1, R2}.

On the contrary, from an analytic point of view, most of the attention
has been given to the operators related to the CR-structure on $H_n$,  
or to its
sub-Riemannian structure (the sub-Laplacian and the Kohn Laplacians),  
leaving only a marginal
r\^ole to the ``Riemannian'' operators.
Our interest here is in the Hodge Laplacians
$\Delta_k=dd^*+d^*d$ acting on differential $k$-forms on $H_n$, a
family of operators that naturally arise in the Riemannian setting, and  
in their
$L^p$-functional calculus.
For $k\ge1$, $\Delta_k$ is far from being
diagonal (in contrast with the Kohn Laplacians for the
$\bar\de_b$-complex) in any reasonable basis of forms. This makes its  
analysis quite involved, with
a level of complexity that increases with $k$ (as long as $k\le n$; it
goes without saying that we are dispensed from treating higher values  
of $k$ by
Hodge duality). For this reason our results are limited to the case  
$k=1$ (together with the
``scalar'' case $k=0$), and we believe that investigating Laplacians on  
higher-order forms
would require a more sophisticated understanding of the decomposition  
of the space of
$L^2$-forms under the action of $\Delta_k$.

Our main result is Theorem 6.8, proving that, if $m$ is a
Mihlin-H\"ormander multiplier on the positive half-line with a  
sufficiently high order of
smoothness, then $m(\Delta_1)$ is bounded on 1-forms in $L^p$, for  
$1<p<\infty$. The
order of smoothness is measured in terms of ``scale-invariant'' local
Sobolev norms (called $L^2_{\tau,\text{\rm sloc}}$), and $\tau$ is  
required to be strictly larger than
$n+\half$, i.e. half of the dimension of $H_n$ as a manifold.

As a preliminary result, the same statement is proved for the
Laplace-Beltrami operator $\Delta_0$ acting on functions (Theorem
6.4). That the critical value for $\tau$ is $n+\half$ in this case is  
not surprising, because
$\Delta_0$ locally behaves like the ordinary Laplacian on $\R^{2n+1}$,  
and at
infinity like the sub-Laplacian $L$, and it is known that $n+\half$ is  
critical
for both these operators (see \cite{MS} for what concerns $L$). To be
more specific, if we scale on $H_n$ isotropically by a parameter  
tending to zero, we
produce a deformation of $\Delta_0$ which in the limit gives the
Laplacian; on the other hand,  if we scale by a parameter tending to  
infinity in the automorphic
(non-isotropic) way, the resulting deformation of $\Delta_0$ tends to
$L$ \cite{NRS}.
As observed in \cite{R2}, this doubly asymptotic picture has no
analogue for forms of order $k\ge1$. The fact that $n+\half$ remains  
the critical value for
$\tau$ also when $k=1$ turns out to be a consequence of the fact that
the space of $L^2$-1-forms  decomposes as the orthogonal sum of five
subspaces such that on each of them the action of $\Delta_1$ is  
unitarily equivalent (possibly modulo an intertwining
operator) to the action of a ``scalar'' differential or
pseudo-differential operator related to $\Delta_0$. Precisely, we find   
\roster \item
the space $V_1$ of exact forms, where the action of $\Delta_1$ is
unitarily equivalent to $\Delta_0$ acting on scalar functions;
\item the space $V_2^+$ of $\de^*_b$-closed $(1,0)$-forms, where
$\Delta_1$ acts as $\Delta_0-iT$ componentwise; \item the space $V_2^-$  
of $\bar\de^*_b$-closed $(0,1)$-forms, where $\Delta_1$ acts as  
$\Delta_0+iT$ componentwise; \item two other subspaces, $V_3^\pm$,   
where the action of $\Delta_1$
is unitarily equivalent to that of $\Delta_0+\frac
n2\pm\sqrt{\Delta_0+\frac{n^2}4}$ on scalar functions. \endroster

Whereas $V_1, V_2^+,V_2^-$ can be detected by a simple inspection, the
last two subspaces are not so visible, and their description involves a  
rather delicate
formalism. The  presence of $V_3^\pm$ had been detected before in
\cite{L} for $H_1$. We thank Michael Christ for bringing this reference  
to our
attention.

Once this is established, the task is to prove first that the
decomposition of the space  of 1-forms into these five subspaces also  
makes sense in $L^p$ for
$p\ne2$ in the range  $1<p<\infty$ (i.e. to prove that the
corresponding orthogonal projections are $L^p$-bounded), and then to
prove that Mihlin-H\"ormander multipliers with order of smoothness  
$\tau>n+\half$ give bounded operators on $L^p$ when applied
to the five operators above. In doing so, we heavily rely on the
results in \cite{MRS1, MRS2}.
\vskip.6cm

\head 1. The Hodge Laplacians\endhead

\vskip.3cm

Let $H_n$ be the $(2n+1)$-dimensional Heisenberg group with
coordinates $(x,y,t)\in \R^n\times\R^n\times\R$, and with a basis of  
left-invariant vector
fields given by $$ X_j=\de_{x_j}-\frac{y_j}2\de_t\ ,\qquad
Y_j=\de_{y_j}+\frac{x_j}2\de_t\ ,\qquad T=\de_t\ ,\tag1.1
$$ for $1\le j\le n$. The dual basis of 1-forms is given by the $2n$
elementary forms $dx_j, \ dy_j$ and by the contact form
$$
\theta=dt-\half\sum_{j=1}^n(x_jdy_j-y_jdx_j)\ .
$$

We denote by $\Lambda^k=\Lambda^k(\h_n^*)$ the $k$-th exterior product
of the dual of the Lie algebra $\h_n$  of $H_n$ (also identifiable
with the space of left-invariant $k$-forms on $H_n$). We call
$$
\D\Lambda^k(H_n)=\D(H_n)\otimes \Lambda^k
$$ the space of smooth $k$-forms on $H_n$ with compact support. This
notation will be consistently adapted to function spaces other than
$\D(H_n)$ or to subspaces of $\Lambda^k$.

We shall often meet differential (or pseudo-differential) operators
which act separately on each scalar component of a given form. In
these cases we will denote by the same symbol the operator, call it  
$D$, acting on scalar-valued functions, and the corresponding operator
acting on forms, which should be more correctly denoted by $D\otimes I$.

It will be convenient for us to work with different bases of complex
vector fields and forms. We then set
$$ B_j=\frac1{\sqrt2}(X_j-iY_j)\ ,\qquad \bar  
B_j=\frac1{\sqrt2}(X_j+iY_j)\ ,
$$ and
$$
\beta_j=\frac1{\sqrt2}(dx_j+idy_j)\ ,\qquad
\bar\beta_j=\frac1{\sqrt2}(dx_j-idy_j)\ . $$

The relevant commutation relation is
$$ [B_j,\bar B_j]=iT\ .
$$

The differential $df$ of a smooth function is then given by
$$ df=\sum_{j=1}^n(X_jf\,dx_j+Y_jf\,dy_j)+Tf\,
\theta=\sum_{j=1}^n(B_jf\,\beta_j+\bar B_jf\,\bar\beta_j)+Tf\,\theta\  
,\tag1.2 $$ and similarly for exterior derivatives of differential  
forms. Observe
that, in particular,
$$ d\theta=-\sum_{j=1}^ndx_j\wedge dy_j=-i\sum_{j=1}^n \beta_j\wedge
\bar\beta_j\ .\tag1.3
$$

A $k$-form  $\omega$ decomposes uniquely as
$$
\omega=\omega_1+\theta\wedge \omega_2\ ,\tag1.4
$$ with
$$
\aligned
\omega_1&=\sum_{|I|+|I'|=k}f_{I,I'}\beta^I\wedge \bar\beta^{I'}\\
\omega_2&=\sum_{|I|+|I'|=k-1}g_{J,J'}\beta^J\wedge \bar\beta^{J'}\ ,
\endaligned\tag1.5
$$ where we have followed the usual convention that, if
$I=\{i_1,\dots,i_p\}$ is a finite subset of $\{1,\dots,n\}$ with  
$i_1<i_2<\cdots<i_p$, then $$
\beta^I=\beta_{i_1}\wedge \beta_{i_2}\wedge\cdots\wedge\beta_{i_p}\ ,
$$ and similarly for $\bar \beta^{I'}$.

Clearly, $\omega_2=0$ if $k=0$ and $\omega_1=0$ for $k=2n+1$.

We say that $\omega$ is {\it horizontal} if $\omega_2=0$, and we call
{\it horizontal differential} of a smooth function $f$ the horizontal
form $$ d_Hf=\sum_{j=1}^n(B_jf\,\beta_j+\bar B_jf\,\bar\beta_j)\  
.\tag1.6
$$

We denote by $\Lambda^k_H$ (resp. $\D\Lambda^k_H(H_n)$) the subspace
of $\Lambda^k$ consisting of horizontal $k$-forms which are  
left-invariant
(resp. with compact support).

The notion of ``horizontal form'' presents serious problems, that are
treated in a systematic way in \cite{R1}. For instance, the natural
extension to horizontal forms of the operator $d_H$ in (1.6) does not  
define a complex, because
$d_H^2\ne0$. However we shall not use any such property, and on the
other hand (1.6) provides a convenient notation. For instance, w.r. to  
the decomposition (1.4), we have    $$
d(\omega_1+\theta\wedge\omega_2)
=\big(d_H\omega_1+(d\theta)\wedge\omega_2  
\big)+\theta\wedge(T\omega_1-d_H\omega_2)\ .\tag1.7 $$

Identifying $\omega$ with the pair $\pmatrix \omega_1\\
\omega_2\endpmatrix$ in (1.4), the oper\-ator $$  
d=d_k:\D\Lambda^k(H_n)\rightarrow \D\Lambda^{k+1}(H_n)
$$ is then repre\-sented by the matrix  $$ d =  \pmatrix  
d_H&e(d\theta)\\ T&-d_H \endpmatrix\ ,\tag1.8
$$ where $e$ denotes exterior multiplication, i.e.  
$e(d\theta)\omega=(d\theta)\wedge\omega$.
\vskip.3cm

We introduce on $H_n$ the left-invariant Riemannian metric that makes
the basis (1.1) orthonormal at each point. W.r. to the induced inner  
product on
$\Lambda^k$, the elements
$$
\beta^I\wedge\bar\beta^{I'}\ ,\qquad
\theta\wedge\beta^J\wedge\bar\beta^{J'} $$ (with $|I|+|I'|=k$,  
$|J|+|J'|=k-1$) also form an orthonormal
basis. Let $\omega,\omega'$ be two $k$-forms, with
$$
\omega=\omega_1+\theta\wedge\omega_2\ ,\qquad
\omega'=\omega'_1+\theta\wedge\omega'_2\ , $$ and let  
$f_{I,I'},g_{J,J'}$ be the coefficients of $\omega$ as in
(1.5),  and $f'_{I,I'},g'_{J,J'}$ the corresponding coefficients of  
$\omega'$. The
inner product in $L^2\Lambda^k(H_n)=L^2(H_n)\otimes \Lambda^k$ is such  
that $$
\align
\lan\omega,\omega'\ran_k&=\sum_{I,I'}\lan f_{I,I'},f'_{I,I'}\ran
+\sum_{J,J'}\lan g_{J,J'},g'_{J,J'}\ran\\
&=\lan\omega_1,\omega'_1\ran_k+\lan\omega_2,\omega'_2\ran_{k-1}\ ,
\endalign
$$ where the inner products of the coefficients are taken in
$L^2(H_n)$. In particular the decomposition (1.4) is orthogonal. The  
formal adjoint of
$d_{k-1}$,  $$ d^*=d_{k-1}^*:\D\Lambda^k(H_n)\rightarrow  
\D\Lambda^{k-1}(H_n)
$$ is represented by the adjoint matrix of (1.8), i.e.
$$ d^* =  \pmatrix d^*_H&-T\\i(d\theta)&-d^*_H \endpmatrix\ ,\tag1.9 $$  
where $i(d\theta)=e(d\theta)^*$ is the interior multiplication
operator $$ i(d\theta)\omega= i\sum_{j=1}^n  
i(\bar\beta_j)i(\beta_j)\omega =
i\sum_{j=1}^n \bar\beta_j\lrcorner(\beta_j\lrcorner\omega)\ .
$$

It follows that the Hodge Laplacian on $k$-forms
$$
\Delta_k=dd^*+d^*d
$$ is expressed by the matrix
$$
\aligned
\Delta_k&=\pmatrix d_H&e(d\theta)\\ T&-d_H \endpmatrix  \pmatrix  
d^*_H&-T\\i(d\theta)&-d^*_H
\endpmatrix +  \pmatrix d^*_H&-T\\i(d\theta)&-d^*_H \endpmatrix
\pmatrix d_H&e(d\theta)\\ T&-d_H \endpmatrix\\ \\ &=\pmatrix  
\Delta_H-T^2+e(d\theta)i(d\theta) &
\big[d_H^*,e(d\theta)\big]\\ \big[i(d\theta),d_H\big]&
\Delta_H-T^2+i(d\theta)e(d\theta)\endpmatrix\ , \endaligned\tag1.10
$$ where
$$
\Delta_H= d_Hd_H^*+d_H^*d_H\ .
$$

In particular, for $k=0$ we simply have
$$
\Delta_0=d^*d=-\sum_{j=1}^n (B_j\bar B_j+\bar B_jB_j)-T^2 \ ,\tag1.11
$$ acting on scalar-valued functions.

\vskip.6cm

\head 2. The CR-structure  \endhead

\vskip.3cm

It is possible to simplify various terms in (1.10) and get a better
understanding of that formula by appropriate decompositions of the  
space of horizontal
forms. In order to do so, we must refer to the standard CR-structure on  
$H_n$. The
operators $$
\de_bf=\sum_{j=1}^n B_jf\,\beta_j\ ,\qquad \bar\de_bf=\sum_{j=1}^n\bar
B_jf\,\bar\beta_j\ ,
$$ initially defined on functions, are naturally extended to forms. They
satisfy the following identities:
$$
\de_b^2=\bar\de_b^2= \de_b\bar\de_b^*+\bar\de_b^*\de_b= \bar \de_b
\de_b^*+ \de_b^*\bar \de_b=0\ ,\tag2.1
$$ as well as
$$ d_H=\de_b+\bar\de_b\ .\tag2.2
$$
Observe that, by (1.8),
$$ d_H^2=\de_b\bar\de_b+\bar\de_b\de_b=-Te(d\theta)\ .
$$

Setting
$$
\square=\de_b\de_b^*+\de_b^*\de_b \ ,\qquad \overline\square  
=\bar\de_b\bar\de_b^*+\bar\de_b^*\bar\de_b \ ,
$$ we obtain that
$$
\Delta_H=\square+\overline\square\ .\tag2.3
$$

Then $\overline\square$ is the Kohn Laplacian and $\square$ its
complex conjugate.
A {\it $(p,q)$-form}, $p,q\le n$, is a horizontal form
$$
\omega=\sum_{|I|=p\,,|I'|=q}f_{I,I'}\beta^I\wedge \bar\beta^{I'}\ .
$$

Clearly, the decomposition (1.4) can be further refined, by
decomposing $\omega_1$ as a sum of $(p,k-p)$-forms and $\omega_2$ as a  
sum of
$(p,k-1-p)$-forms. The notation $\Lambda^{p,q}$,  
$\D\Lambda^{p,q}(H_n)$, etc. refers to
$(p,q)$-forms.
It is well known \cite{FS} that $\square$ and $\overline\square$ act
as scalar operators on
$(p,q)$-forms (we shall write $\square_{p,q}$ and
$\overline\square_{p,q}$ when appropriate). If
$$ L=-\sum_{j=1}^n (B_j\bar B_j+\bar B_jB_j)=-\sum_{j=1}^n(X_j^2+Y_j^2)
$$ is the sub-Laplacian, then
$$
\square_{p,q}=\half L+i\Big(\frac n2-p\Big)T\ ;\tag2.4
$$ similarly,
$$
\overline\square_{p,q}=\half L-i\Big(\frac n2-q\Big)T\ .\tag2.5
$$
It follows from (2.3) that
$$
\Delta_H=L+i(q-p)T\tag2.6
$$ on $(p,q)$-forms.

\vskip.3cm

We next describe the structure of the remaining diagonal terms in
(1.10), i.e. $e(d\theta)i(d\theta)$ and its transpose
$i(d\theta)e(d\theta)$. Since these operators do not involve any
differentiation, their action can be analyzed on exterior forms. Many
of the formulas below are also stated in \cite{R1,2} and are derived  
from the formulas for the Lefshetz
decomposition on K\"ahler manifolds in
\cite{W}. For completeness, we give some explicit proofs, and we allow
forms of any order, even though we shall later restrict ourselves to
1-forms.
\proclaim{Proposition 2.1} Consider the following subspaces of  
$\Lambda^{p,q}$, $$
\align V^{p,q}_j&=e(d\theta)^j\ker_{\Lambda^{p-j,q-j}}i(d\theta)\ ,\\
W^{p,q}_\ell&=i(d\theta)^\ell\ker_{\Lambda^{p+\ell,q+\ell}}e(d\theta)\ .
\endalign
$$
Then $V^{p,q}_j$ is non-trivial if and only if $\max\{0,k-n\}\le
j\le\min\{p,q\}$, $W^{p,q}_\ell$ is non-trivial if and only if  
$\max\{0,n-k\}\le
\ell\le\min\{n-p,n-q\}$, and we have the equality
$$ V^{p,q}_j=W^{p,q}_\ell\ ,\qquad \text{ for\ } \ell=j+n-k=\ell(j)\ .
$$
Moreover, $\Lambda^{p,q}$  is the orthogonal sum of the non-trivial
$V^{p,q}_j$, and $$
\align
e(d\theta)i(d\theta)
&=j(j+1+n-k)=\big(\ell(j)+1\big)\big(\ell(j)+k-n\big)\\  
i(d\theta)e(d\theta)
&=(j+1)(j+n-k)=\ell(j)\big(\ell(j)+1+k-n\big)
\endalign
$$ on $V^{p,q}_j$. \endproclaim

\proof Because $\ker
i(d\theta)=\big(e(d\theta)\Lambda^{p-1,q-1}\big)^\perp$ inside  
$\Lambda^{p,q}$, every $(p,q)$-form $\om$ can be uniquely decomposed
into the orthogonal sum $$
\om=\om_0+e(d\theta) \alpha
$$ with $i(d\theta)\om_0=0$. Next, we decompose $\al$ as
$$
\alpha=\om_1+e(d\theta)\alpha'\ ,
$$ with $i(d\theta)\om_1=0$. The resulting decomposition
$$
\om=\om_0+e(d\theta)\om_1+e(d\theta)^2\alpha'
$$ is also orthogonal, because
$$
\align
\lan e(d\theta)\om_1,e(d\theta)^2\alpha'\ran&=
\lan i(d\theta)e(d\theta)\om_1,e(d\theta)\alpha'\ran\\
&=(n-k+2)\lan\om_1,e(d\theta)\alpha'\ran\\ &=0\ .
\endalign
$$

Iterating this procedure, we end up with writing $$
\om=\sum_{j=0}^{\min\{p,q\}}e(d\theta)^j\om_j\ ,
$$ with $\om_j\in\Lambda^{p-j,q-j}$ and $i(d\theta)\om_j=0$.
A direct computation shows that, when applied to horizontal $k$-forms,  
$$
\big[i(d\theta),e(d\theta)]=(n-k)I\ ,\tag2.7
$$ (see also \cite{W}). Then
$$
\aligned i(d\theta)e(d\theta)^j\om_j
&=[i(d\theta),e(d\theta)^j]\om_j\\
&= \sum_{i=0}^{j-1}e(d\theta)^i[i(d\theta),
e(d\theta)]e(d\theta)^{j-1-i}\om_j\\  
&=\sum_{i=0}^{j-1}(n-k+2+2i)e(d\theta)^{j-1}\om_j\\   
&=j(n-k+j+1)e(d\theta)^{j-1}\om_j\
.
\endaligned\tag2.8
$$

Hence,
$$ e(d\theta)i(d\theta)\om
=\sum_{j=0}^{\min\{p,q\}} j(n-k+j+1) e(d\theta)^j\om_j\
,\tag2.9
$$ showing that $ e(d\theta)i(d\theta)$ diagonalizes w.r. to the
decomposition $$
\Lambda^{p,q}=\sum_{j=0}^{\min\{p,q\}} V^{p,q}_j\ .
$$

By (2.7), $i(d\theta)e(d\theta)$ also diagonalizes w.r. to the same
decomposition, and $$  \aligned i(d\theta)e(d\theta)\om
&=e(d\theta)i(d\theta)\om+(n-k)\om\\
&=\sum_{j=0}^{\min\{p,q\}} (j+1)(n-k+j) e(d\theta)^j\om_j\ .
\endaligned\tag2.10
$$

But $i(d\theta)e(d\theta)$ is positive semidefinite, so that $\om_j$
must be 0 for $j<k-n$. Therefore $V^{p,q}_j$ can be non-trivial only if  
$\max\{0,k-n\}\le
j\le\min\{p,q\}$. In order to see that this condition is also  
sufficient, observe that for $j$ in this range, $0\le p+q-2j\le
k-2\max\{0,k-n\}=\min\{k,2n-k\}\le n$.
Then $$
\om=\beta_1\wedge\cdots\wedge
\beta_{p-j}\wedge\bar\beta_{p-j+1}\wedge\cdots
\wedge\bar\beta_{p+q-2j} $$
is a non-zero element of $\Lambda^{p-j,q-j}$ satisfying
$i(d\theta)\om=0$. That $e(d\theta)^j\om$ is non-zero is trivial for  
$j=0$ and it follows by induction from
(2.8). In conclusion,  $$
\Lambda^{p,q}=\sum_{j=\max\{0,k-n\}}^{\min\{p,q\}} V^{p,q}_j\ ,
$$
where the summands are non-trivial and mutually orthogonal.

A repetition of the same arguments with the r\^oles of $e(d\theta)$
and $i(d\theta)$ interchanged shows that
$$
\Lambda^{p,q}=\sum_{\ell=\max\{0,n-k\}}^{\min\{n-p,n-q\}} W^{p,q}_\ell\  
,
$$
and that
$i(d\theta)e(d\theta)=\ell(\ell+k-n+1)I$ on $W^{p,q}_\ell$.

A comparison with the eigenvalues in (2.10) provides the
identification of $V^{p,q}_j$ with
$W^{p,q}_{\ell(j)}$.
\endproof

\vskip.3cm
Consider now the off-diagonal terms
$$
\big[i(d\theta),d_H\big]\ ,\qquad \big[d_H^*,e(d\theta)\big]
$$
in (1.10). They can be simplified using the following identities.

\proclaim{Proposition 2.2} We have
$$
\gather
\big[i(d\theta),\de_b\big]
= -i\bar\de_b^*\ ,\qquad \big[i(d\theta),\bar\de_b\big]=i\de_b^*\\
\big[\de_b^*,e(d\theta)\big]
= i\bar\de_b\ ,\qquad \big[\bar\de_b^*,e(d\theta)\big]=-i\de_b \ .
\endgather
$$
In particular,
$$
\big[i(d\theta),d_H\big]
=i\de_b^*-i\bar\de_b^*\ ,\qquad \big[d_H^*,e(d\theta)\big]
=i\bar\de_b-i\de_b\ .
$$
\endproclaim

\proof  Given $j\in\{1,\dots,n\}$ and $I,J\subseteq\{1,\dots,n\}$,  
define $\eps_{j,I}^J$ as 0 unless $j\not\in I$ and $\{j\}\cup I= J$, in  
which
case $$
\eps_{j,I}^J=\prod_{i\in I}\sgn(i-j)\ ,
$$
i.e. the signature of the permutation that moves $j$ from the left of
$I$ to its correct position w.r. to the natural ordering of $J$.
Let $\omega=f\,\beta^I\wedge\bar\beta^{I'}$. Then
$$
\align
\big(i(\beta_j)\de_b &+\de_bi(\beta_j)\big)\om
=i(\beta_j)\sum_{\ell,J} \eps_{\ell,I}^J B_\ell
f\, \beta^J\wedge\bar\beta^{I'}\\ &\qquad +\de_b \sum_M  
\eps_{j,M}^If\,\beta^M\wedge\bar \beta^{I'}\\ & =\sum_{\ell,J,L}  
\eps_{\ell,I}^J\eps_{j,L}^J B_\ell
f\,\beta^L\wedge\bar\beta^{I'}\\  &\qquad +\sum_{\ell,M,L}  
\eps_{j,M}^I\eps_{\ell,M}^L
B_\ell f\,\beta^L\wedge\bar \beta^{I'}\\ &=\sum_{\ell  
,L}\bigg(\sum_J\eps_{\ell
,I}^J\eps_{j,L}^J+ \sum_M \eps_{j ,M}^I\eps_{\ell ,M}^L\bigg) B_\ell
f\,\beta^L\wedge\bar \beta^{I'}\ . \endalign $$

Consider the expression
$$
\sum_J\eps_{\ell ,I}^J\eps_{j,L}^J+ \sum_M \eps_{j,M}^I\eps_{\ell ,M}^L
$$
for fixed $\ell, L$. Assume first that $\ell\ne j$. The first sum does
not vanish only in one case: $\ell\not\in I$, $j\in I$,  
$L=I\cup\{\ell\}\setminus\{j\}$,
with the only non-vanishing term in the sum corresponding to  
$J=\{\ell\}\cup I$. But
this is also the only case when the second sum has a non-vanishing  
term, the one
corresponding to $M=I\cap L$. It takes a few moments to verify that, if  
this is the case, the two terms
have opposite signs, so that the total expression is always 0 for  
$\ell\ne j$.

Assume now that $\ell=j\in I$. The first term is 0, and the second
term is also 0 unless $M=I\setminus\{j\}$  and $L=I$. In this case the  
total expression
gives 1.
Finally, if $\ell=j\not\in I$, the first term is 1 and the second is
0. The conclusion is that the expression under consideration equals 1  
if $\ell=j$ and $L=I$
and 0 otherwise. Hence $$ \big(i(\beta_j)\de_b+\de_bi(\beta_j)\big)\om  
=B_j f\,\beta^I\wedge\bar
\beta^{I'}\ . $$

A similar computation shows that
$$
\big(i(\bar\beta_j)\de_b+\de_bi(\bar\beta_j)\big)\om =0\ .
$$

Putting these identities together, we find that
$$
\align
i(d\theta)\de_b\om
&=i\sum_{j=1}^n i(\bar\beta_j)i(\beta_j)\de_b\om\\
&= i\sum_{j=1}^n B_jf i(\bar\beta_j) \beta^I\wedge\bar \beta^{I'} -
i\sum_{j=1}^n i(\bar\beta_j) \de_b i(\beta_j) \om\\
&= i\sum_{j=1}^n B_jf i(\bar\beta_j) \beta^I\wedge\bar \beta^{I'} +
i\sum_{j=1}^n \de_b i(\bar\beta_j) i(\beta_j) \om\\
&= -i\bar\de_b^*\om +\de_b i(d\theta)\om\ .
\endalign
$$

This gives the first identity in the statement. Taking complex
conjugates and transposes, the other three follow. \endproof

In combination with the formula preceding (2.3), this immediately gives 

\proclaim{Corollary 2.3} We have 
$$\align
\square \bar\de_b&=\bar\de_b\square -iT\bar\de_b,\\
\overline\square \de_b&=\de_b\overline\square +iT\de_b,
\endalign
$$
hence, by duality,
$$\align
\bar\de_b^*\square&=\square\bar\de_b^* -iT\bar\de_b^*,\\
\de_b^*\overline\square&=\overline\square\de_b^* +iT\de_b^*.
\endalign
$$
\endproclaim 
\vskip.6cm

\head 3. Spectral multipliers of $i\inv T$ and $L$ \endhead

\vskip.3cm

The operators $i\inv T$ and $L$ admit commuting self-adjoint
extensions on $L^2(H_n)$, and their joint spectrum is the {\it
Heisenberg fan} $F_n\subset\R^2$. If $$
\ell_m=\{(\la,\xi):\xi=(n+2m)|\la| ,\la\in\R\}\ ,
$$
then
$$
F_n=\overline{\bigcup_{m\in\N}\ell_m}\ .
$$

The variable $\la$ corresponds to $i\inv T$ and $\xi$ to $L$, i.e.,
calling $dE(\la,\xi)$ the spectral measure on $F_n$, $$
i\inv T=\int_{F_n}\la\,dE(\la,\xi)\ ,\qquad
L=\int_{F_n}\xi\,dE(\la,\xi)\ . $$

It follows from the Plancherel formula that the spectral measure of
the vertical half-line $\{(0,\xi):\xi\ge0\}\subset F_n$ is zero. A  
spectral multiplier is therefore a function $\mu(\la,\xi)$ on $F_n$  
whose restriction to each $\ell_m$ is measurable w.r. to $d\la$.
Later on we shall use results from \cite{MRS1,2} concerning  
$L^p$-boundedness of
spectral multipliers. For the moment, we use these facts to discuss  
$L^2$-boundedness of certain
operators that will appear in the next Section, together with some
$L^p-L^q$-estimates for unbounded multipliers.

\proclaim{Lemma 3.1} The operators
$$
L^r(\Delta_0+i\al T)^{-r}\ , \qquad T^{2r}(\Delta_0+i\al T)^{-r}
$$
are bounded on $L^2(H_n)$ for $|\al|<n$ and $r>0$.
\endproclaim

\proof By (1.11), $\Delta_0=L-T^2$. Hence we just need to observe that
the multipliers $$
\mu_1(\la,\xi)
=\frac {\xi^r}{(\xi+\la^2-\al\la)^r}\ ,\qquad \mu_2(\la,\xi)=\frac
{\la^{2r}}{(\xi+\la^2-\al\la)^r} $$
are bounded on $F_n$.
\endproof

The {\it Cauchy-Szeg\"o projection} $C$ is the orthogonal projection
of $L^2(H_n)$ onto the Hardy space $H^2(H_n)$, consisting of the
$L^2$-functions $f$ such that $\bar\de_bf=0$. It is a well-known fact
(see \cite{S, Ch. XIII}) that $H^2(H_n)$ is also the null-space of $$  
L-inT=2\overline\square_{0,0}=2\bar\de_b^*\bar\de_b =-2\sum_{j=1}^n  
B_j\bar B_j\ .
$$

What is relevant for us at this stage is that $C=\mu(i\inv T,L)$,
where $\mu$ is equal to 1 on the half-line $\xi=-n\la$, with $\la<0$,
and 0 elsewhere.
In the same way, the complex conjugate $\bar C$ of $C$ projects
$L^2(H_n)$ onto the null space of $\de_b$, which is the same as the  
null space of
$L+inT=2\square_{0,0}$, and its multiplier equals 1 on the half-line  
$\xi=n\la$, with $\la>0$, and 0
elsewhere. The next result follows easily.

\proclaim{Lemma 3.2} The operators
$$
L^r(L-inT)^{-r}(I-C) \ , \qquad L^r(L+inT)^{-r}(I-\bar C)
$$
are bounded on $L^2$.
\endproclaim

We pass now to the $L^p-L^q$-estimates.

\proclaim{Lemma 3.3} Let $\mu(\la,\xi)$ be a smooth function defined
on an angle $D_\del=\{(\la,\xi):\xi>(n-\del)|\la|\}$, with $\del>0$,
and homogeneous of degree $-d$, with $0<d<n+1$. Then $\mu(i\inv T,L)$
is well-defined and bounded from $L^p(H_n)$ to $L^q(H_n)$ for
$1<p<q<\infty$ and $\frac1p-\frac1q=\frac d{n+1}$. \endproclaim

\proof It follows from \cite{G, AD} that $\mu(i\inv T,L)f=f*K$, where
$K$ is smooth  away from the origin and homogeneous of degree  
$-(2n+2-2d)$. The
conclusion follows from the generalized Young inequality.
\endproof

\vskip.6cm

\head 4. Decomposition of $L^2\Lambda^1(H_n)$ under the action of
$\Delta_1$  \endhead
For $k=1$, the conclusions of Sections 1 and 2 lead us to write the
generic 1-form $\om$ as $$
\om=\om_++\om_-+h\theta\ ,
$$
where $\om_+$ is a $(1,0)$-form and $\om_-$ is a $(0,1)$-form.  Then $$
\Delta_1\pmatrix \om_+\\ \om_-\\ h\endpmatrix=
\pmatrix \Delta_0-iT&0&-i\de_b\\
0&\Delta_0+iT&i\bar\de_b\\
i\de_b^*& -i\bar\de_b^* &\Delta_0+n\endpmatrix
\pmatrix \om_+\\ \om_-\\ h\endpmatrix\ .\tag4.1
$$

Obviously, $\Delta_1$, initially defined on $\D\Lambda^1(H_n)$, is
essentially self-adjoint, and the domain of its self-adjoint
(Friedrichs) extension is $$
\dom\Delta_1=\big\{\om\in L^2\Lambda^1(H_n):\Delta_1\om\in
L^2\Lambda^1(H_n)\big\}\ ,
$$
where $\Delta_1\om$ is meant in the sense of distributions.

If $\om\in\D\Lambda^1(H_n)$ is exact, say $\om=d\ph$, then
$$
\Delta_1\om=dd^*d\ph=d\Delta_0\ph\ ,
$$
i.e. $d$ intertwines the action of $\Delta_1$ on $\om$ with the action
of $\Delta_0$ on $\ph$. We shall show that a similar statement holds  
for exact $L^2$-forms, with
$d$ replaced by a modified intertwining operator which is
$L^2$-bounded. Before doing so,  we must make some preliminary
remarks.
\proclaim{Lemma 4.1} The operator  $R=d\Delta_0^{-\half}$
is isometric from $L^2(H_n)$ to its image in $L^2\Lambda^1(H_n)$.   
\endproclaim

\proof By Lemma 3.1, $L^\half\Delta_0^{-\half}$ and
  $T\Delta_0^{-\half}$ are bounded on $L^2(H_n)$.  We recall that the  
Riesz transforms $B_jL^{-\half}$, $\bar
  B_jL^{-\half}$ are also bounded on $L^2$. Since $$
B_j\Delta_0^{-\half}=(B_jL^{-\half})(L^\half\Delta_0^{-\half})\ ,
$$
it follows that $B_j\Delta_0^{-\half}$ is $L^2$-bounded, and similarly
for $\bar B_j\Delta_0\inv$.

Hence, for $\ph\in\D(H_n)$,
$$
\align
\|R\ph\|_2^2&=\|d\Delta_0^{-\half}\ph\|_2^2\\
&=\lan\Delta_0^{-\half} d^*d\Delta_0^{-\half}\ph,\ph\ran\\
&=\|\ph\|_2^2\ .\qquad\qed
\endalign
$$
\enddemo

We say that $\om\in L^2\Lambda^1(H_n)$ is {\it exact} if there exists  
$u\in\D'(H_n)$ such that $\om=du$ in the sense of distributions
(componentwise). In the same sense we shall talk later on of {\it  
$\de_b$-exact} forms or of
{\it $\bar\de_b$-exact} forms.
\proclaim{Lemma 4.2} Let $r$ be such that
$\half-\frac1r=\frac1{2n+2}$. If $\om\in L^2\Lambda^1(H_n)$ is exact,  
then $\om=dv$, in the sense of
distributions, for some $v\in L^r(H_n)$.
\endproclaim

\proof By definition, there is $u\in\D'(H_n)$ such that
$\om=du$. Define $$
v=\Delta_0\inv d^*\om=L^{-\half}(\Delta_0^{-\half} L^\half)R^*\om\ ,
$$
where $R^*=\Delta_0^{-\half}d^*$ is the adjoint of the operator $R$ in
Lemma 4.1. Then $R^*$ is $L^2$-bounded, and $\Delta_0^{-\half}
L^\half$ is too, by the spectral theorem. Finally, $L^{-\half}$ is  
bounded from $L^2$ to $L^r$, e.g. by Lemma
3.3. Hence $v\in L^r(H_n)$. Moreover,
$$
dv=RR^*\om\in L^2\Lambda^1(H_n)\ ,
$$
and
$$
\Delta_0v=d^*\om=\Delta_0u\ .
$$

Observe now that
$$
\Delta_1(\om-dv)=\Delta_1d(u-v)=d\Delta_0(u-v)=0\ .
$$

The conclusion will follow from the next lemma.
\endproof

\proclaim{Lemma 4.3} The Hodge Laplacian $\Delta_1$ is injective on  
$L^2\Lambda^1(H_n)$.
\endproclaim

\proof Assume that $\om=\om_++\om_-+h\theta$ satisfies $\Delta_1\om=0$
in the sense of distributions. By (4.1), this means that
$$
\aligned
(\Delta_0-iT)\om_+&=i\de_bh\ ,\\
(\Delta_0+iT)\om_-&=-i\bar\de_bh\ ,\\
(\Delta_0+n)h&=-i\de_b^*\om_++i\bar\de_b^*\om_-\ .
\endaligned\tag4.2
$$

We multiply the first equation in (4.2) by $(\Delta_0-iT)\de_b^*$, and
the second equation by $(\Delta_0+iT)\bar\de_b^*$. Using the
identities  $$
\de_b^*(\Delta_0-iT)=(\Delta_0+iT)\de_b^*\ ,\qquad  
\bar\de_b^*(\Delta_0+iT)=(\Delta_0-iT)\bar\de_b^*
$$
-- easily deduced from (2.4)
and (2.5) --, and performing some simple computations, we obtain that
$$
(\Delta_0^2+T^2)(\Delta_0+n)h=\big(\Delta_0^2+T^2(\Delta_0+n)\big)h\ ,
$$
i.e. $$
\Delta_0^2(\Delta_0+n-1)h=0\ .
$$

Since the zero set of the multiplier corresponding to the operator on
the left-hand side is the origin, and it has  measure zero in the  
Heisenberg fan, this
implies that $h=0$. \endproof

\proclaim{Proposition 4.4} The operator $P_1=RR^*$ on  
$L^2\Lambda^1(H_n)$  is the
orthogonal projection onto the subspace of exact $L^2$-forms. In
particular,  this  subspace is closed. Moreover, $P_1$ maps  
$\dom\Delta_1$ into itself. \endproclaim

\proof Clearly, $P_1$ is self-adjoint. Assume that $\om\in  
L^2\Lambda^1(H_n)$ is exact. By Lemma 4.2, there
is $v\in L^r(H_n)$ such that $\om=dv$.
Let $\chi$ be a non-negative, smooth function on $H_n$ with compact
support, equal to 1 on a neighborhood of the origin, and define
$\chi_j(z,t)=\chi(z/j,t/j^2)$. Let also $\{\ph_j\}_{j\in\N}$ be an  
approximate identity in $\D(H_n)$. If  $$ v_j=\ph_j*(\chi_jv)\ ,
$$
then $v_j\rightarrow v$ in $L^r(H_n)$. Moreover,
$$
dv_j=\ph_j*(\chi_j\om)+\ph_j*(vd\chi_j)\ ,
$$
if we interpret the concolution $\ph_j*\al$ of $\ph_j$ with a 1-form  
$\al$
componentwise.

If $|\cdot |$ denotes a homogeneous norm on $H_n$,
$$
\align
\|\ph_j*(vd\chi_j)\|_2&\le\|vd\chi_j\|_2\\
&\le \frac Cj \bigg(\int_{|x|\sim j}|v|^2\bigg)^\half\\
&\le \frac Cj (j^{2n+2})^{\frac{r-2}{2r}}\bigg(
\int_{|x|\sim j}|v|^r\bigg)^{\frac1r}\\ &=C\bigg(\int_{|x|\sim  
j}|v|^r\bigg)^{\frac1r}\ ,
\endalign
$$
and it tends to zero as $j$ tends to infinity. Hence $dv_j\rightarrow
dv$ in $L^2\Lambda^1(H_n)$.

Given $\sigma\in\D\Lambda^1(H_n)$, we then have
$$
\align
\lan P_1\om,\sigma\ran&=\lan dv,P_1\sigma\ran\\
&=\lim_{j\to\infty}\lan dv_j,P_1\sigma\ran\\
&=\lim_{j\to\infty}\lan v_j, d^*d\Delta_0\inv d^*\sigma\ran\\
&=\lim_{j\to\infty}\lan v_j,  d^*\sigma\ran\\
&=\lan v,d^*\sigma\ran\\
&=\lan dv,\sigma\ran\ ,
\endalign
$$
showing that $P_1\om=\om$.

On the other hand, if $\om=P_1\om'$, let $v=\Delta_0\inv d^*\om'\in
L^r(H_n)$, as in the proof of Lemma 4.2. If
$\sigma\in\D\Lambda^1(H_n)$, $$
\align
\lan dv,\sigma\ran&=\lan \Delta_0\inv d^*\om',d^*\sigma\ran\\
&=\lan\om',P_1\sigma\ran\ ,
\endalign
$$
so that $dv=P_1\om'=\om$.
To prove the last part of the statement, take again
$\sigma\in\D\Lambda^1(H_n)$. Then $$
\Delta_1P_1\sigma=\Delta_1 (d\Delta_0\inv d^*\sigma)=dd^*d\Delta_0\inv
d^*\sigma=dd^*\sigma\ ,
$$
and
$$
P_1\Delta_1\sigma=d\Delta_0\inv d^*\Delta_1\sigma=d\Delta_0\inv
d^*dd^*\sigma=dd^*\sigma\ .
$$

Therefore $\Delta_1P_1=P_1\Delta_1$ on $\D\Lambda^1(H_n)$. For a
general $\sigma\in\dom\Delta_1$, we take a sequence of forms
$\sigma_j\in\D\Lambda^1(H_n)$ such that
$\sigma_j\to\sigma$ and $\Delta_1\sigma_j\to\Delta_1\sigma$ in the
$L^2$-norm. Then $P_1\sigma_j\to P_1\sigma$, and
$$
P_1\Delta_1\sigma=\lim_{j\to\infty}P_1\Delta_1\sigma_j=\lim_{j\to\infty} 
\Delta_1P_1
\sigma_j\ .
$$

Since $\Delta_1$ is closed, $P_1\sigma\in\dom\Delta_1$, and
$P_1\Delta_1\sigma=\Delta_1P_1\sigma$.
\endproof

\proclaim{Proposition 4.5} Let $V_1$ be the range of $P_1$ in  
$L^2\Lambda^1(H_n)$, i.e.
the space of exact $L^2$-forms. Then $R$ maps $\dom\Delta_0$ onto
$(\dom\Delta_1)\cap V_1$, and intertwines the action of $\Delta_0$ with  
that of
$\Delta_1$, i.e.
$$
R\Delta_0=\Delta_1R\ ,
$$
on $\dom\Delta_0$.
\endproclaim

\proof If $\ph\in\D(H_n)$,
$$
\align
\Delta_1R\ph&=(dd^*+d^*d)d\Delta_0^{-\half}\ph\\
&=d(d^*d)\Delta_0^{-\half}\ph\\
&=d\Delta_0^{-\half}\Delta_0\ph\\
&=R\Delta_0\ph\ .
\endalign
$$

An adaptation of the proof of Proposition 4.4 shows that  
$R(\dom\Delta_0)\subseteq\dom\Delta_1$, and that $R\Delta_0=\Delta_1R$  
on $\dom\Delta_0$. 

Conversely, take $\om\in(\dom\Delta_1)\cap V_1$ and
$\ph\in\D(H_n)$. Since $\D\Lambda^1(H_n)$ is a core for $\Delta_1,$ we
find a sequence $\{\om_j\}_j$ in this space such that $\om=\lim \om_j$
and $\Delta_1\om=\lim \Delta_1\om_j$ in $L^2.$ Moreover, $\Delta_1
R\ph=R\Delta_0\ph\in L^2\Lambda^1(H_n),$ and thus 
$$
\align
\lan R^*\Delta_1\om,\ph\ran&=\lan\Delta_1\om,R\ph\ran\\
&=\lim_{j\to\infty} \lan\Delta_1\om_j,R\ph\ran=\lim_{j\to\infty}
\lan\om_j,\Delta_1 R\ph\ran\\
 &=\lim_{j\to\infty}\lan\om_j, R\Delta_0\ph\ran\\
 &=\lan\om, R\Delta_0\ph\ran\\
&=\lan \Delta_0 R^*\om,\ph\ran\ ,
\endalign
$$
showing that $\Delta_0(R^*\om)$, defined in the sense of
distributions, is equal to $R^*\Delta_1\om$. In particular,  
$R^*\om\in\dom\Delta_0$.

Since $\om=P_1\om=R(R^*\om)$, it follows that $\om\in R(\dom\Delta_0)$.  
\endproof

\vskip.3cm

We are so led to study $\Delta_1$ on $V_1^\perp$, the orthogonal
complement of the exact $L^2$-forms.  This is the space of co-closed  
$L^2$-forms, i.e. the forms $\om$ such
that $d^*\om=0$.
We denote by $V_2^+$ (resp. $V_2^-$) the space of co-closed $(1,0)$
forms (resp. $(0,1)$ forms).

\proclaim{Proposition 4.6} For $\om\in V_2^\pm$,
$\Delta_1\om=(\Delta_0\mp iT)\om$ in the sense of distributions.
\endproclaim

\proof If $\om\in V_2^+$, then $\de_b^*\om=d^*\om=0$. The conclusion
follows from (4.1), and similarly for $V_2^-$. \endproof

Observe that, on $H_1$, $V_2^+$ consists of the $(1,0)$-forms
$f\beta$ with
$\de_b^*(f\beta)=-\bar Bf=0$. Therefore, on $V_2^+$,
$$
\Delta_0-iT=-(2 B\bar B+T^2)=-T^2\ .
$$

In the same way, $\Delta_0+iT=-T^2$ on $V_2^-$.

\vskip.3cm

We want to describe now the orthogonal projections $P_2^\pm$ from
$L^2\Lambda^1(H_n)$ onto $V_2^\pm$. We look at $P_2^+$ as the
composition of the orthogonal projection $Q^+$ from
$L^2\Lambda^1(H_n)$ onto $L^2\Lambda^{1,0}(H_n)$ followed by the
orthogonal projection from $L^2\Lambda^{1,0}(H_n)$ onto $V_2^+$ (and
similarly for $P_2^-$).
It is immediate to verify that
$$
Q^\pm(\om_++\om_-+h\theta)=\om_\pm\ .
$$

In order to describe the second factor in the decomposition of
$P_2^+$, it is preferable to consider its complementary projection,  
from $L^2\Lambda^{1,0}(H_n)$
onto the orthogonal complement $(V_2^+)^\perp$. Since $V_2^+$ is the  
null space of
$\de_b^*$,  $(V_2^+)^\perp$ is the closure of the space of
$\de_b$-exact $L^2$-$(1,0)$-forms. In the same way, $(V_2^-)^\perp$ is  
the closure
of the space of $\bar\de_b$-exact $L^2$-$(0,1)$-forms.

These projections involve the operators
$\square_{0,0},\overline\square_{0,0}$ in (2.4) and (2.5),
$$
\square_{0,0} =\de_b^*\de_b=\half(L+inT)\ ,\qquad \overline\square_{0,0}
=\bar\de_b^*\bar\de_b=\half(L-inT)\ .
$$

As there will be no confusion from now on, we drop the double
subscript and simply write $\square$ and $\overline\square$. As we have  
observed already,
$$
\square u=0
\Longleftrightarrow\de_bu=0\ ,\qquad \overline\square
u=0\Longleftrightarrow\bar\de_bu=0\ .\tag4.3
$$
Consequently, the image in $L^2(H_n)$ of
$\de_b^*$ is contained in
$\big(\ker\square\big)^\perp$ and the image in $L^2(H_n)$ of  
$\bar\de_b^*$ is
contained in
$\big(\ker\overline\square\big)^\perp$.

In particular, $$
\de_b^*=(I-\bar C)\de_b^*\ ,\qquad \bar\de_b^*=(I- C)\bar\de_b^*\ .  
\tag4.4
$$

It follows from Lemma 3.2 and boundedness of the Riesz transforms that  
$$
\align
\square^{-\half}\de_b^*&=\sqrt 2 \Big(L^\half(L+inT)^{-\half}(I-\bar
C)\Big)\,\big(L^{-\half}\de_b^*\big)\ ,\\
\overline\square^{-\half}\bar\de_b^*&=\sqrt 2\Big(L^\half(L-inT)^{- 
\half}(I-
C)\Big)\,\big(L^{-\half}\bar\de_b^*\big)
\endalign
$$
are well defined and bounded from $L^2\Lambda^{1,0}(H_n)$ (resp.
$L^2\Lambda^{0,1}(H_n)$) to $L^2(H_n)$. If the factors $I-C$ and
$I-\bar C$ are superfluous in the above formulas because of (4.4), the
same is not true for the adjoint operators,  
$\de_b\square^{-\half}(I-\bar C)$ and
$\bar\de_b\overline\square^{-\half}(I-C)$.

We conclude that the four operators we will be dealing with,
$$
\aligned
\Ri=\de_b\square^{-\half}(I-\bar C)\ ,
&\qquad
\bar\Ri=\bar\de_b\overline\square^{-\half}(I-C)\ ,\\
\Ri^*=\square^{-\half}\de_b^*\ ,
&\qquad \bar\Ri^*=\overline\square^{-\half}\bar\de_b^*\
,
\endaligned
\tag4.5
$$
are $L^2$-bounded.

\proclaim{Proposition 4.7} The operator $\Ri\Ri^*$ is the orthogonal  
projection from
$L^2\Lambda^{1,0}(H_n)$ onto the subspace of $\de_b$-exact forms, and  
$\bar\Ri\bar\Ri^*$ is the orthogonal projection from
$L^2\Lambda^{0,1}(H_n)$ onto the subspace of $\bar\de_b$-exact forms.  
In particular, these two
subspaces are closed. Moreover, $\Ri^*\Ri=I-\bar C,\ 
\bar\Ri^*\bar\Ri=I-C.$
\endproclaim

\proof The argument is the same as in the proof of Proposition 4.4.   
The only
substantial  difference is that we must write $$
\square\inv\de_b^*=\square^{-\half}(I-\bar C) \square^{-\half}\de_b^*\
, $$
and notice that Lemma 3.3 can be applied to the factor
$\square^{-\half}(I-\bar C)$. In fact this operator can be realized as
$\mu(i\inv T,L)$, if $\mu$ is an appropriately chosen smooth function  
on some angle $D_\del$, homogeneous of degree
$-1/2$, equal to $(\xi-n\la)^{-\half}$ on $F_n$ except for the  
half-line $\xi=n\la$,
$\la>0$, where it is set equal to 0.
\endproof

\proclaim{Corollary 4.8} The orthogonal projections $P_2^\pm$ from  
$L^2\Lambda^1(H_n)$ onto
$V_2^\pm$ are given by $$
P_2^+=(I-\Ri\Ri^*)Q^+\ ,\qquad P_2^-=(I-\bar\Ri\bar\Ri^*)Q^-\ .
$$

They map $\dom\Delta_1$ into itself.
\endproclaim

\vskip.6cm

\head 5. Decomposition of the action of $\Delta_1$ on $V_3$ \endhead

\vskip.3cm

It remains to describe the action of $\Delta_1$ on the orthogonal
complement  $V_3$ of $V_1\oplus V_2^+\oplus V_2^-$ in  
$L^2\Lambda^1(H_n)$. Notice that
$\Delta_1(V_3\cap \dom \Delta_1)\subset V_3$. It follows from
Proposition 4.4 and Corollary 4.8 that  $V_3\cap\dom \Delta_1$ is dense  
in $V_3$.

In order to describe $V_3$ we take a detour that has the advantage of
making this space somewhat more explicit.  We forget for a moment that  
$V_1$
has been analyzed already, and we look at the full orthogonal complement
of $V_2$,
$$
V_2^\perp =
\big\{  \om\,:\, \om=\om_+ +\om_- +h\theta\, , \om_+\,  \text{is}\ 
\p_b-\text{exact}\, ,\  \om_-\,\text{is}\  \bar\p_b-\text{exact}\big\}\ .
$$
Since $\om_+$ is $\p_b$-exact, let $u=\Ri^* \om_+$.  Then $u\in L^2$ and
$\bar Cu=0$.   Moreover, we can recover $\om_+$ from $u$, since
$\om_+=\Ri u$, by Prop. 4.7.  Analogously, we set $v=\bar\Ri^* \om_-$
so that $v\in L^2$, $Cv=0$ and $\om_-=\bar\Ri v$.
Thus, we are lead to consider the closed subspace of
$\big(L^2\bigr)^3$
$$
W = \big\{ (u,v,h)\in \big(L^2\bigr)^3\, :\ \bar Cu=Cv=0\big\} \, .
$$
\proclaim{Lemma 5.1} Define $\Gamma: W  \rightarrow V_2^\perp$ by  
setting
$$
\Gamma (u,v,h) =\Ri u + \bar\Ri v+h\theta\ .
$$
Then $\Gamma$ is unitary and its inverse is given by
$$
\Gamma^*(\om) = \bigl(\Ri^*\om_+,\bar\Ri^*\om_-,h\bigr)\, .
$$
\endproclaim
\proof
By definition of $\Ri$ and $\bar\Ri$ it is clear that $\Gamma$ maps
$W$ into $V_2^\perp$.  Next, by Proposition 4.7, 
$$
\aligned
\lan \Gamma(u,v,h),\Gamma(u',v',h')\ran
& = \lan \Ri u,\Ri u'\ran +
\lan \bar\Ri v,\bar\Ri v'\ran + \lan h,h'\ran \\
& = \lan \Ri^* \Ri u, u'\ran +
\lan \bar\Ri^* \bar\Ri v,v'\ran + \lan h,h'\ran \\
& = \lan(u,v,h),(u',v',h')\ran\ ,
\endaligned
$$
which shows that $\Gamma$ preserves the inner product.
The previous discussion shows that $\Gamma^*
\Gamma=\text{Id}_W$, and furthermore $\Gamma$ is onto since
  $\Gamma
\Gamma^*=\text{Id}_{V_2^\perp}$.
\endproof

We now set
$$
D_1 = \Gamma^* \Delta_1 \Gamma\, ,
$$
being $\dom(D_1)= \Gamma^*\big( \dom(\Delta_1)\cap V_2^\perp\big)$.  We
compute
$D_1$ explicitely.  Writing $\Gamma(u,v,h)=\om(u,v,h)$ and
  recalling that
$\Delta_1$ is given by (4.1), we have
$$
\Delta_1 \om(u,v,h)= \om(u',v',h')\, ,
$$
where
$$
\left\{
\aligned
\Ri u' &= (\Delta_0 -iT)\Ri u -i\p_b h\\
\bar\Ri v' &= (\Delta_0 +iT)\bar\Ri v +i\bar\p_b h\\
h' & = i\p_b^* \Ri u-i\bar\p_b^* \bar\Ri v +(\Delta_0 +n)h\, .
\endaligned
\right.
$$
By applying $\Ri^*$ to the first equation and $\bar\Ri^*$ to the
second one and using the commutation relations from Corollary 2.3  we
obtain 
$$
\left\{
\aligned
u' &= (\Delta_0 +iT)u -i\square^\half h\\
v' &= (\Delta_0 -iT)v +i\overline\square^{\half}  h\\
h' & = i\square^\half u-i\overline\square^\half v +(\Delta_0 +n)h\, .
\endaligned
\right.
$$

Therefore,
$$
D_1= \pmatrix \Delta_0+iT&0&-i\square^\half\\
         0&\Delta_0-iT&i\overline\square^\half\\
         i\square^\half& -i\overline\square^\half &\Delta_0+n          
\endpmatrix \, .
$$
Consider the corresponding matrix of spectral multipliers
$$
\aligned
d_1
& =   \pmatrix \xi+\la^2 -\la &0&-i\sqrt{\thalf(\xi-n\la)}\\
         0&\xi+\la^2+\la&i\sqrt{\thalf(\xi+n\la)}\\
         i\sqrt{\thalf(\xi-n\la)}& -i\sqrt{\thalf(\xi+n\la)}
&\xi+\la^2+n          \endpmatrix \\
& = (\xi+\la^2)I + \pmatrix -\la &0&-i\sqrt{\thalf(\xi-n\la)}\\
         0& \la&i\sqrt{\thalf(\xi+n\la)}\\
         i\sqrt{\thalf(\xi-n\la)}& -i\sqrt{\thalf(\xi+n\la)} & n          
   \endpmatrix \, .
\endaligned
$$
Diagonalization of $d_1$ will have the following implication.
Assume that $$
v=\pmatrix a_1(\la,\xi) \\
a_2(\la,\xi) \\
a_3(\la,\xi) \endpmatrix
$$
is a unit eigenvector of $d_1$ of eigenvalue $\mu(\la,\xi)$. If we take  
a
scalar function $f\in L^2(H_n)$ such that
$$ F=\pmatrix a_1(i^{-1}T,L)f \\
a_2(i^{-1}T,L)f \\
a_3(i^{-1}T,L)f \endpmatrix \in W\ ,
$$
then $D_1F=\mu(i^{-1}T,L)F$.

\proclaim{Lemma 5.2}
The eigenvalues of $d_1$ are $\xi+\la^2$ and $\xi+\la^2+\frac{n}{2}\pm  
\sqrt{\xi+\la^2+\frac{n^2}{4}}$. The matrix entries of the  orthogonal
projections to the  eigenspaces of $d_1$ are functions of
$(\la,\xi)$ which are bounded on the Heisenberg fan $F_n$.
\endproclaim

\proof
We compute the eigenvalues of $m_1=d_1-(\xi+\la^2)I$:
$$
\det(m_1-\mu I) = -\mu^3 +n\mu^2 +(\xi+\la^2)\mu \, ,
$$
so that $m_1$ has eigenvalues
$$
\mu=0\, ,\qquad \mu_\pm=\frac{n}{2}\pm\sqrt{\xi+\la^2+\frac{n^2}{4}}  
\,
. $$

Next we determine the eigenvectors and the orthogonal projections onto
the eigenspaces of $m_1$.

A unit eigenvector for $\mu=0$ is
$$
v_0 = \frac1{\sqrt{\xi+\la^2}}\pmatrix \sqrt{\thalf(\xi-n\la)} \\
         \sqrt{\thalf(\xi+n\la)}\\
         i\la \endpmatrix\ .
$$

In order to describe the eigenvectors  corresponding to the
eigengvalues $\mu_\pm$, we set $$
\aligned
a&=a(\la,\xi)= \sqrt{\xi+\la^2+\frac{n^2}{4}}\\
q_\delta^\eps&=q_\delta^\eps(\la,\xi) =a+\eps\frac{n}{2}+\delta\la\, ,  
\endaligned\tag5.2
$$
where $\eps,\delta=\pm1$.  Notice that the following identities hold:
$$
\aligned
& q_+^+ q_-^- = \xi-n\la \\
& q_-^+ q_+^- = \xi+n\la\\
& q_+^+ + q_-^- = q^+_- + q_+^- = 2a \\
& q^+_+q^-_+=(a+\la)^2-\frac{n^2}4\\
& q_-^+q_-^-=(a-\la)^2-\frac{n^2}4\\
& q_+^+q_-^+=(a+\frac n2)^2-\la^2\\
& q_+^-q_-^-=(a-\frac n2)^2-\la^2.
\endaligned
\tag5.3
$$

Since
$$
\aligned
& m_1-\mu_\pm I= \\
& \pmatrix -\la-\frac{n}{2}\mp\sqrt{\xi+\la^2+\frac{n^2}{4}}
& 0 & -i\sqrt{\thalf(\xi-n\la)} \\
0 & \la -\frac{n}{2}\mp\sqrt{\xi+\la^2+\frac{n^2}{4}}
&   i\sqrt{\thalf(\xi+n\la)} \\
i\sqrt{\thalf(\xi-n\la)} & -i\sqrt{\thalf(\xi+n\la)} & \frac{n}{2}
  \mp\sqrt{\xi+\la^2+\frac{n^2}{4}}
\endpmatrix \, ,
\endaligned
$$
a unit eigenvector relative to $\mu_+$ is
$$
v_+  
=\frac1{\sqrt {2a(a+\frac n2)}} \pmatrix -i\sqrt{\half
q^+_-q^-_-}
\\ i\sqrt{\half q^+_+q^-_+} \\
\sqrt{q^+_+q^+_-} \endpmatrix\ ,
$$
where we have used the identities (5.3) to obtain the normalizing  
factor.

Similar computations show that
a unit eigenvector relative to $\mu_-$ is
$$
v_- =\frac1{\sqrt {2a(a-\frac n2)}} \pmatrix i\sqrt{\half  
q^+_+q^-_+}
\\ -i\sqrt{\half q^+_-q^-_-} \\
\sqrt{q^-_+q^-_-} \endpmatrix
\ .
$$

The  orthogonal
projections corresponding to the eigenvectors above are represented by
the matrices
$p_0=v_0v_0^*$, $p_\pm=v_\pm v^*_\pm$. Clearly these three matrices
satisfy
$p_0+p_++p_-=I$, and their  entries
are bounded by 1 on the fan $F_n$. \endproof

Next, we wish to decompose $W$ as the direct sum of  subspaces  
in such a way that $D_1$ acts as a scalar operator on any of these
subspaces.
Recalling the definition (5.2) of $a$ and
$q_\delta^\eps$, we set
$$
\Cal A=\sqrt{\Delta_0+\frac{n^2}{4}}=a(i^{-1}T,L) \ ,\qquad Q_\delta^\eps
= \Cal A +\eps\frac{n}{2}-\delta iT=q^\eps_\delta(i^{-1}T,L)\, , \tag5.6
$$ where $\eps,\delta=\pm1$.  By (5.3) and (2.4), (2.5) we then have the
following identities
$$
\aligned
& Q_+^+ Q_-^- = 2\square \\
& Q_+^- Q_-^+ = 2\overline\square\ .\
 \endaligned
\tag5.7
$$

\proclaim{Proposition 5.3}
Define $S_0$, $S_\pm$ resp. to be the operators from $L^2(H_n)$ to
${L^2(H_n)}^3$ having
$v_0$, $v_\pm$ resp. as spectral multipliers. Then $S_0$ and $S_+$ map
$L^2(H_n)$ isometrically into $W$, and $S_-$ maps $L^2_0(H_n)=\{f\in
L^2(H_n):Cf=\bar Cf=0\}$ isometrically into $W$.

Moreover, $W$ is the orthogonal sum of $W_0=S_0L^2(H_n)$,  
$W_+=S_+L^2(H_n)$,
$W_-=S_-L^2_0(H_n)$. More precisely, every $(u,v,h)\in W$ decomposes
uniquely as $$
(u,v,h)=S_0f_0+S_+f_++S_-f_-\ ,
$$
with $f_0=S_0^*(u,v,h)\in L^2(H_n)$, $f_+=S_+^*(u,v,h)\in L^2(H_n)$, and
$f_-=S_-^*(u,v,h)\in L^2_0(H_n)$.
Finally, the operators $P_0=S_0S_0^*$, $P_\pm=S_\pm S_\pm^*$ on
$W$ whose spectral multipliers are
$p_0$, $p_\pm$ resp., are the orthogonal projections onto $W_0$, $W_\pm$
respectively.
\endproclaim
\proof
We know that, for every fixed $(\la,\xi)\in F_n,$ 
$$I=p_0+p_++p_-=v_0v_0^*+v_+v_+^*+v_-v_-^*\quad\text{on}\ \C^3,
$$
where $p_0,p_+,p_-$ are pairwise orthogonal projections. By the spectral
theorem, this implies 
$$I=P_0+P_+ +P_-=S_0S_0^*+S_+S_+^*+S_-S_-^*\quad\text{on}\ L^2(H_n)^3,
$$
where $P_0,P_+$ and $P_-$ are pairwise orthogonal projections. Moreover,
since the spectral multiplier for $S_0^*S_0$ is $v_0^*v_0=||v_0||^2=1,$ 
$S_0$ is isometric, and the same is true for
$S_+,S_-,$  by similar reasoning.

Thus, every $(u,v,h)\in L^2(H_n)^3$ uniquely decomposes as the orthogonal
sum 
$$(u,v,h)=S_0f_0+S_+f_++S_-f_-,$$
with $f_0=S_0^*(u,v,h), f_+=S_+^*(u,v,h), f_-=S_-^*(u,v,h)\in L^2(H_n)$.
There remains to prove that the mapping 
$$\Cal T:(f_0, f_+,f_-)\mapsto S_0f_0+S_+f_++S_-f_-,$$
when restricted to the subspace $\Omega=L^2(H_n)\times L^2(H_n)\times
L_0^2(H_n),$ maps into and onto $W.$
To this end, notice that the first components in $v_0$ and $v_+$ vanish
for  
$\xi=n\la$.
Together with the fact that the spectral multiplier of $\bar C$ is
the characteristic function of the set where $\xi=n\la$, this implies  
that,
if
$(u,v,h)$ equals $S_0f_0$ or $S_+f_+$, then
$\bar Cu=0$. A similar argument shows that $Cv=0$. The
same conclusion holds for $(u,v,h)=S_-f_-$ if we impose that
$Cf_-=\bar Cf_-=0$, i.e. $f_-\in L^2_0(H_n)$. Thus $\Cal T(\Omega)\subset
W.$

Conversely, given $(u,v,h)\in W$, define $f_0,f_\pm\in L^2(H_n)$ as in
the   statement.
In particular,
$$
f_-=\frac1{\sqrt{2\Cal A\big(\Cal A -\frac n2\big)}} \bigg(
i\sqrt{\half Q^+_+Q^-_+}u -i\sqrt{\half
Q^+_-Q^-_-}v+\sqrt{Q^-_+Q^-_-}h\bigg)\ .
$$
 From the identities
$$
\bar Cu=0\ ,\qquad Cv=0\ ,\qquad \bar CQ^-_-=0\ ,\qquad CQ^-_+=0\ ,
$$
we conclude that $\bar Cf_-=Cf_-=0,$ hence $(f_0,f_+,f_-)\in\Omega.$
\endproof

\proclaim{Remark}
{\rm It can be proved that it is possible to give another description of
the three subspaces of $W$ as
$$ \aligned
& W_0  =   \bigl\{ (u,v,h)\in W\, : Tu=\square^\half h,\  
Tv=\overline\square^\half h \bigr\} \\
& W_+ =   \bigl\{ (u,v,h)\in W\, : Q_+^+ u =-i\square^\half h,\ Q_-^+  
v=i \overline\square^\half h \bigr\}  \\
& W_- =   \bigl\{ (u,v,h)\in W\, : Q_-^- u=i\square^\half h,\ Q_+^-  
v=-i\overline\square^\half h \bigr\} \ . \endaligned
$$
\endproclaim

\vskip.5cm

Composing with $\Gamma$, this decomposition of $W$ gives rise to an
orthogonal decomposition of $V_2^\perp$. Notice that, if $(u,v,h)\in
W_0$, i.e.
$$
u=\Delta_0^{-\half}\square^\half f\ ,\qquad
v=\Delta_0^{-\half}\overline\square^\half f\ ,\qquad
h=T\Delta_0^{-\half}f\ ,
$$
for some $f\in L^2(H_n)$, then
$$
\Gamma (u,v,h)=\Ri\Delta_0^{-\half}\square^\half f +\bar\Ri
\Delta_0^{-\half}\overline\square^\half f+T\Delta_0^{-\half}f\theta=Rf\ ,
$$
so that $\Gamma W_0=V_1$, the space of exact forms.

Define
$$
V_3^\pm = \Gamma(W_\pm)\ .
$$
\proclaim{Proposition 5.4}
The orthogonal complement of $V_1\oplus V_2$ in $L^2\Lambda_1$ is the
subspace $V_3=V_3^+\oplus V_3^-$.  The operators $\Gamma S_+$ and
$\Gamma S_-$ are unitary respectively from $L^2(H_n)$ onto $V_3^+$ and
from $L^2_0(H_n)$ onto $V_3^-$.
The orthogonal projections from $V_2^\perp$ onto $V_3^+$ and $V_3^-$
are
$$
\Pi_\pm=\Gamma S_\pm S_\pm^*\Gamma^*\ .
$$

Moreover,
$$
\Gamma S_\pm\Big(\Delta_0+\frac
n2\pm\textstyle{\sqrt{\Delta_0+\frac{n^2}4}}\,\Big)=\Delta_1 \Gamma
S_\pm\ .
$$
\endproclaim

\proof
The first part of the statement is obvious.
What concerns the action of $\Delta_1$ follows from the fact that,
since
$D_1$ commutes with
$P_0$,
$P_\pm$,
$$
\dom (D_1) = \big(\dom (D_1)\cap W_0\big) \oplus
\big(\dom (D_1)\cap W_+ \big)\oplus
\big(\dom (D_1)\cap W_- \big)\ . \quad\qed
$$
\enddemo

\vskip.6cm

\head 6. $L^p$-boundedness of spectral multipliers of $\Delta_1$  
\endhead

\vskip.3cm

On the basis of the previous analysis, we can say that
$$
\Delta_1= R\Delta_0R \inv=R\Delta_0R^*
$$
  on $V_1$,
$$
\Delta_1=\Delta_0\mp iT
$$
on $V_2^\pm$, and, by Proposition 5.4,
$$
\Delta_1=\Gamma S_\pm\Big(\Delta_0+\frac n2\pm
\sqrt{\Delta_0+\frac{n^2}4}\Big)S^*_\pm \Gamma^*
$$
on $V_3^\pm$.

This implies that, given a bounded Borel function $m$ on
$(0,+\infty)=\R^*_+$, the operator $T_m=m(\Delta_1)$ equals $$
Rm(\Delta_0)R^*\ ,\qquad m(\Delta_0\mp iT)\ ,\qquad
\Gamma S_\pm m\Big(\Delta_0+\frac n2\pm
\sqrt{\Delta_0+\frac{n^2}4}\Big)S^*_\pm
\Gamma^* \ ,
$$
on the corresponding subspaces. Denoting by $P_3=I-P_1-P_2^+-P_2^-$
the orthogonal projection from $L^2\Lambda^1(H_n)$ onto $V_3$, we  
obviously have
$$
m(\Delta_1)
=m(\Delta_1)P_1+m(\Delta_1)P_2^++m(\Delta_1)P_2^-+m(\Delta_1)\Pi_+P_3
+m(\Delta_1)\Pi_-P_3\
.
$$

Observe that, since $R^*R=I$, we have $R^*P_1=R^*RR^*=R^*$. Similarly,
$$
S^*_\pm \Gamma^*\Pi_\pm=S^*_\pm \Gamma^*\ .
$$

We then have
$$
\aligned
m(\Delta_1)
&= Rm(\Delta_0)R^*+m(\Delta_0-iT)P_2^++m(\Delta_0+iT)P_2^- \\
&\quad
+\Gamma S_+
m\Big(\Delta_0+\frac n2+ \sqrt{\Delta_0+\frac{n^2}4}\Big)S^*_+ \Gamma^*
P_3\\ &\quad+\Gamma S_-
m\Big(\Delta_0+\frac n2-\sqrt{\Delta_0+\frac{n^2}4}\Big)
S^*_-\Gamma^* P_3\ .
\endaligned\tag6.1
$$

We are so led to discuss $L^p$ boundedness of each of the operators
appearing in (6.1). For this purpose, we recall the following result,  
taken from
\cite{MRS2, Cor.2.4}, and concerning Marcinkiewicz multipliers of  
$i\inv T$ and $L$. We shall
present a series of technical lemmas in a rather self-contained  
fashion. We do not claim
full originality for every single statement. In particular, various  
overlappings with
arguments in \cite{MS} are present.
Given $\rho,\sigma>0$, we say that a function $f(\la,\xi)$ is in the
mixed  Sobolev space $L^2_{\rho,\sigma }=L^2_{\rho,\sigma }(\R^2)$ if
$$
\aligned
\|f\|_{L^2_{\rho,\sigma
}}^2:&=\int_{\R^2}(1+|\xi'|)^{2\rho}(1+|\la'|+|\xi'|)^{2\sigma}|\hat
f(\la',\xi')|^2\,d\la'\,d\xi'\\  
&=c\|(1+|\de_\xi|)^\rho(1+|\de_\la|+|\de_\xi|)^\sigma
f\|_2^2 \endaligned\tag6.2 $$
is finite. When $\rho$ and $\sigma$ are integers,  this condition means  
that the
derivatives $\de^i_\la\de^j_\xi f$ are in $L^2$ for $i\le \sigma$ and  
$i+j\le
\rho+\sigma$. We shall make use of this characterization, together with  
the fact that
the $L^2_{\rho,\sigma }$ form an interpolation family.

Let $\eta\in \D\big((\R^*_+)^2\big)$ be a non-trivial, non-negative,  
smooth function (briefly, a bump function). We say that a bounded  
function $\mu(\la,\xi)$ defined
on $(\R^*_+)^2$ is in $L^2_{\rho,\sigma,\text{\rm  
sloc}}\big((\R^*_+)^2\big)$ if for every
$r=(r_1,r_2)\in(\R^*_+)^2$, the function
$\mu^r(\la,\xi)=\mu(r_1\la,r_2\xi)\eta\big(\la,\xi\big)$ is in  
$L^2_{\rho,\sigma }$
and $$ \|\mu\|_{L^2_{\rho,\sigma,\text{\rm sloc}}}=\sup_r\|\mu^r\|_{
L^2_{\rho,\sigma}}\tag6.3  $$ is finite. We extend this definition to  
functions
$\mu$ defined on $\R\times\R^*_+$ by requiring that both $\mu(\la,\xi)$  
and
$\mu(-\la,\xi)$ are in $L^2_{\rho,\sigma,\text{\rm  
sloc}}\big((\R^*_+)^2\big)$.

If $\rho$ and $\sigma$ are integers, to require that $\mu\in  
L^2_{\rho,\sigma,\text{\rm sloc}}\big((\R^*_+)^2\big)$ is equivalent to  
requiring that
$$
\sup_{r_1,r_2>0}r_1^{-1+2i}r_2^{-1+2j}
\int_{r_1<|\la|<2r_1\,,\,r_2<\xi<2r_2}\big|\de_\la^i\de_\xi^j
\mu(\la,\xi)\big|^2\,d\la\,d\xi <+\infty\ ,\tag6.4 $$ for all $i,j$  
such that $i\le \sigma$, $i+j\le \rho+\sigma$. In
particular, the definition of $L^2_{\rho,\sigma,\text{\rm sloc}}$ is  
independent of
the choice of $\eta$. The same is true for every $\rho,\sigma>0$, as  
the following
lemma shows.

\proclaim{Lemma 6.1} Given two bump functions $\eta_1$ and $\eta_2$ on  
$(\R^*_+)^2$, the
norms $(6.3)$ that they define are equivalent.
Let $\Omega$ be a family of bump functions, such that all
the $\eta\in\Omega$ are supported on the same compact subset of  
$(\R^*_+)^2$, and
that their $C^k$-norms are uniformly bounded for some $k\ge  
\rho+\sigma$. Given another
bump function $\eta_0$, the norms $(6.3)$ defined by the  
$\eta\in\Omega$  are controlled
uniformly by the norm defined by $\eta_0$. \endproclaim

\proof If $\ph\in\D(\R^2)$, the operation of multiplication by $\ph$ is  
continuous on $L^2_{\rho,\sigma}$, with a norm controlled by the  
$C^k$-norm of $\ph$, if $k\ge\rho+\sigma$. This is trivial if $\rho$  
and $\sigma$ are integers, and it follows by interpolation in the  
general case.

Given $\eta_1$ and $\eta_2$ as above, there are $r\rinv1,\dots,r\rinv  
k\in (\R^*_+)^2$
such that $$
\psi(\la,\xi)=\sum_{j=1}^k \eta_1 (r\rinv j_1\la,r\rinv j_2\xi)\ge  
\delta >0
$$
on the support of $\eta_2$. Hence $\eta_2=\ph\psi$ for some  
$\ph\in\D(\R^2)$. Then
$$
\aligned
\|\mu(r\cdot)\eta_2\|_{ L^2_{\rho,\sigma}}&\le C\|\mu(r\cdot)\psi\|_{  
L^2_{\rho,\sigma}}\\
&\le C\sum_{j=1}^k \|\mu(r\cdot)\eta_1(r\rinv j\cdot)\|_{  
L^2_{\rho,\sigma}}\\
&\le C'\sum_{j=1}^k \big\|\mu\big(r(r\rinv  
j)\inv\cdot\big)\eta_1\big\|_{
L^2_{\rho,\sigma}}\ , \endaligned\tag6.5
$$
and this implies the first part of the statement.

Given a family $\Omega$ of bump functions as above, the same $\psi$ can  
be used for all
the $\eta\in\Omega$, because of the condition on the supports. It  
follows that the set
$\{\ph=\eta/\psi:\eta\in\Omega\}$ is bounded in $C^k$ for every $k$.  
Hence the constant
$C'$ appearing in (6.5), with $\eta_2=\eta$ and $\eta_1=\eta_0$, can be  
taken
independently of $\eta$.  \endproof
\proclaim{Theorem 6.2 \cite{MRS2}} Let $\mu$ be a bounded function in
$L^2_{\rho,\sigma,\text{\rm sloc} }(\R\times\R^*_+)$ for some $\rho>n$  
and
$\sigma>\half$. Then $\mu(i\inv T,L)$ is bounded on $L^p(H_n)$ for  
$1<p<\infty$,
with norms controlled by $\|\mu\|_{L^2_{\rho,\sigma,\text{\rm sloc}}}$.  
\endproclaim

 From this statement we shall derive a result concerning spectral  
 multipliers of
$\Delta_0+i\al T$ for $|\al|<n$. Observe that, if $m$ is a bounded  
function on $\R^*_+$,
then $m(\Delta_0+i\al T)=\mu(i\inv T,L)$, with $$  
\mu(\la,\xi)=m(\la^2+\xi-\al\la)\ .
$$

If $\tau>0$, we say that $m\in L^2_{\tau,\text{\rm sloc}}(\R^*_+)$ (or  
that $m$ is a {\it
Mihlin-H\"ormander multiplier of order} $\tau$) if $$
\|m\|_{ L^2_{\tau,\text{\rm  
sloc}}}=\sup_{r>0}\big\|m(r\cdot)\ph\big\|_{L^2_\tau}
$$
is finite, where $\ph$ is a non-trivial, non-negative, smooth bump  
function on $\R^*_+$ and the $L^2_\tau$-norm is the ordinary Sobolev  
norm on $\R$.
It will be useful to observe that $m\in L^2_{\tau,\text{\rm  
sloc}}(\R^*_+)$ if and only
if $\mu(\la,\xi)=m(\xi)$ is in $L^2_{\tau,\sigma,\text{\rm sloc}  
}\big((\R^*_+)^2\big)$
for any $\sigma$. In particular, the analogue of Lemma 6.1 can be  
formulated, with the
obvious modifications, for $L^2_{\tau,\text{\rm sloc}}(\R^*_+)$.

One important technical aspect of our argument is the following.

\proclaim{Proposition 6.3} Let $\rho,\sigma>0$, $\al\in(-n,n)$, and let  
$m$ be a
Mihlin-H\"ormander multiplier of order $\tau=\rho+\sigma$. Then
$\mu(\la,\xi)=m(\la^2+\xi-\al\la)$ coincides on $F_n$ with a function in
$L^2_{\rho,\sigma,\text{\rm loc}}$. \endproclaim

This and Theorem 6.2 imply the following result.

\proclaim{Theorem 6.4} If $m$ is a Mihlin-H\"ormander multiplier of  
order
$\tau>n+\half$, then $m(\Delta_0+i\al T)$ is a bounded operator on  
$L^p(H_n)$ for
$|\al|<n$ and
$1<p<\infty$.
\endproclaim

The proof of Proposition 6.3 requires a few lemmas.

\proclaim{Lemma 6.5} If $\mu(\la,\xi)$ is in  
$L^2_{\rho,\sigma,\text{\rm sloc}}\big((\R^*_+)^2\big)$, then  
$\mu(\la^2,\xi)$ is in $L^2_{\rho,\sigma,\text{\rm  
sloc}}(\R\times\R^*_+)$, with equivalence of norms.
\endproclaim
\proof Let $K$ be a compact subset of $(\R^*_+)^2$.  If $\rho$ and  
$\sigma$ are
integers, it is quite clear that a function $f(\la,\xi)$ supported on  
$K$ is in
$L^2_{\rho,\sigma }$ if and only if $f(\la^2,\xi)$ is in  
$L^2_{\rho,\sigma}$. By complex
interpolation, the same holds for all $\rho,\sigma>0$.

In order to prove the Lemma, it is sufficient to consider the  
restriction
$\tilde\mu(\la,\xi)$ of $\mu(\la^2,\xi)$ to $(\R^*_+)^2$. If $\eta$ is  
a bump function,
the $L^2_{\rho,\sigma}$-norm of  $$
\tilde\mu(r_1\la,r_2\xi)\eta(\la,\xi)=  
\mu(r_1^2\la^2,r_2\xi)\eta(\la,\xi)
$$
is controlled by the $L^2_{\rho,\sigma}$-norm of
$\mu(r_1^2\la,r_2\xi)\eta(\sqrt\la,\xi)$. The conclusion follows easily  
from Lemma 6.1.
\endproof

\proclaim{Lemma 6.6}  Let $\mu\in L^2_{\rho,\sigma,\text{\rm  
sloc}}\big((\R^*_+)^2\big)$,
and $\del>0$. Let also $\psi$ be smooth on $\R\times\R^*_+$,
homogeneous of degree zero, and supported on the angle
$D_\delta=\{(\la,\xi):\xi\ge(n-\delta)|\la|\}$. If $\al<n-\delta$, then  
$$
\mu'(\la,\xi)=\mu(\la,\xi-\al\la)\psi(\la,\xi) $$ is also in  
$L^2_{\rho,\sigma,\text{\rm
sloc}}\big((\R^*_+)^2\big)$. \endproclaim

\proof If $\gamma\in\R$, the linear change of variables  
$(\la,\xi)\longmapsto
(\la,\xi+\gamma\la)$ induces an isomorphism of $L^2  
_{\rho,\sigma}(\R^2)$ onto itself,
with constants controlled by $\gamma$. This follows easily from (6.2).
Therefore, if $\eta_0$ is a bump function on $(\R^*_+)^2$ and  
$r_1,r_2>0$, the
$L^2_{\rho,\sigma}$-norm of $$
\mu'(r_1\la,r_2\xi)\eta(\la,\xi)=\mu(r_1\la,r_2\xi-\al
r_1\la)\psi(r_1\la,r_2\xi)\eta_0(\la,\xi) $$
is equivalent to the $L^2_{\rho,\sigma}$-norm of
$$
\mu(r_1\la,r_2\xi')\psi(r_1\la,r_2\xi'+\al
r_1\la)\eta_0\Big(\la,\xi'+\al\frac{r_1}{r_2}\la\Big)\ , $$
with constants controlled by the ratio $r_1/r_2$. Let
$$
\eta_r(\la,\xi')= \psi(r_1\la,r_2\xi'+\al
r_1\la)\eta_0\Big(\la,\xi'+\al\frac{r_1}{r_2}\la\Big)=
\psi\Big(\la,\frac{r_2}{r_1}\xi'+\al
\la\Big)\eta_0\Big(\la,\xi'+\al\frac{r_1}{r_2}\la\Big)\ . $$

The conclusion follows if we prove that, for an appropriate choice of  
$\eta_0$, the set
$\Omega=\{\eta_r:r\in(\R^*_+)^2\}$ satisfies the assumptions of Lemma  
6.1, and that
$\eta_r\ne0$ only if the ratio $r_1/r_2$ is bounded.

Assume that the support of $\eta_0$ is contained in the square  
$[1-\eps,1+\eps]^2$,
with $\eps\in(0,1)$ to be determined. A necessary condition for having  
$\eta_r\ne0$
is that there exists $(\la,\xi')$ such that the conditions
$$
\Big(\la,\xi'+\al\frac{r_1}{r_2}\la\Big)\in[1-\eps,1+\eps]^2\ ,\qquad
  \Big(\la,\frac{r_2}{r_1}\xi'+\al \la)\in D_\delta\ ,\tag6.6
$$
are satisfied, or, otherwise stated, that
$$
\bigg([1-\eps,1+\eps]\times\Big[\frac{r_2}{r_1}(1-\eps),  
\frac{r_2}{r_1}(1+\eps)\Big]\bigg) \cap D_\delta\ne\emptyset\ .
$$

This occurs if and only if the point $\Big(1-\eps,  
\frac{r_2}{r_1}(1+\eps)\Big)\in D_\delta$, i.e. if and only if
$$
\frac{r_2}{r_1}>(n-\delta)\frac{1-\eps}{1+\eps}\ .\tag6.7
$$

The requirement about the boundedness of the ratios $r_1/r_2$ is then  
fulfilled.

Once (6.7) is satisfied, we check that the supports of the $\eta_r$   
are contained in a
common compact subset of $(\R^*_+)^2$. Clearly, if  
$(\la,\xi')\in\supp\eta_r$, then
$\la\in[1-\eps,1+\eps]$. As to $\xi'$, we impose, for all the  
$(r_1,r_2)$ satisfying
(6.7), the condition $$ 1-\eps\le \xi'+\al\frac{r_1}{r_2}\la\le 1+\eps\  
, $$
taken from the first of (6.6). The existence of an upper bound for  
$\xi'$ independent of
$r$ follows from the fact that the ratios $r_1/r_2$ are bounded. For  
the lower bound,
there is no problem if $\al\le0$. If $0<\al<n-\delta$, taking into  
account (6.7) and
that $\la<1+\eps$, we are done if   $$
(1-\eps)-\al\frac{1+\eps}{(n-\delta)(1-\eps)}(1+\eps)=
(1-\eps)\bigg(1-\frac\al{n-\delta}\Big(\frac{1+\eps}{1- 
\eps}\Big)^2\bigg)
> 0\ .
$$

This can be obtained by choosing $\eps$ small enough. A simple  
verification shows that the derivatives of the $\eta_r$ are uniformly  
bounded,
so that the conclusion follows from Lemma 6.1. \endproof
We can now prove Proposition 6.3.

\demo{Proof of Proposition 6.3} Given a Mihlin-H\"ormander multiplier  
$m$ of order
$\tau$ on the positive half-line, consider
$$
\mu_1(\la,\xi)=m(\la+\xi)
$$
on $(\R^*_+)^2$. Applying Lemma 6.6 to $\mu_0(\la,\xi)=m(\xi)$, with
$\del=n$ and $\al=-1$, we obtain that
$\mu_1\in L^2_{\rho,\sigma,\text{\rm sloc}}\big((\R^*_+)^2\big)$.

By Lemma 6.5,
$\mu_2(\la,\xi)=m(\la^2+\xi)$ is in $L^2_{\rho,\sigma,\text{\rm  
sloc}}(\R\times\R^*_+)$.

Given $\al\in\R$ with $|\al|<n$, take $\delta>0$, $\delta<n-|\al|$ and
construct  $\psi$ smooth, homogeneous of degree 0 supported on  
$D_\delta$ and equal to 1
on $F_n$. Applying Lemma 6.6 to $m_2(\pm\la,\xi)$ restricted to  
$(R^*_+)^2$, we conclude that also
$\mu_3(\la,\xi)=m(\la^2+\xi-\al\la)\psi(\la,\xi)$ is in  
$L^2_{\rho,\sigma,\text{\rm
sloc}}(\R\times\R^*_+)$. \endproof

\proclaim{Proposition 6.7} The operators $S_+$ and $S_-$ are bounded from
$L^p(H_n)$ to $L^p(H_n)^3$ for $1<p<\infty.$ 
\endproclaim
\proof The components of the operators $S_+$ and $S_-$ are
spectral multiplier operators whose multipliers are the components of
$v_+$ and $v_-$.

It turns out that the half-lines $\xi=n\la$ and $\xi=-n\la$ play a
special role here, which is why we discuss them separately. We restrict
ourselves to the half-line $\xi=n\la;$ the other half-line can be treated
in a similar way. 

On the former half-line, we have $a=\la+\frac n2,$ and, using (5.3), one
finds that 
$$
v_+  = \pmatrix  0 \\
 i\sqrt{\frac{\la}{2\la+n}} \\
\sqrt{\frac{n}{\la+n}} \endpmatrix \ ,\quad 
v_-=\pmatrix  i\\
0 \\
0 
\endpmatrix \ .
$$
The components of $v_+$ and $v_-$ are Mihlin-H\"ormander multipliers as
functions of $\la>0,$ and since the operator $\bar C,$ which corresponds
to the restriction to the spectral half-line $\xi=n\la,$ is $L^p$-bounded,
we see that the components of $v_+,v_-,$ when restricted to this
half-line, give rise to $L^p$-bounded operators for $1<p<\infty.$ 

In view of the definition of the Heiseberg fan $F_n,$ it thus suffices
to consider the domain where $\xi>(n+1)|\la|.$ Notice that the components
of $v_\pm$ are all  products of multipliers of
the form
$$
\nu_0=(q^\eps_\del)^\half a^{-\half}\
,\qquad 
\nu_+=(q^\eps_\del)^\half \big(a+\textstyle{\frac  
n2}\big)^{-\half}\
,\qquad 
\nu_-=(q^-_\del)^\half \big(a-\textstyle{\frac  
n2}\big)^{-\half}\ .
$$
We show that for $\xi>(n+1)|\la|$ they satisfy the pointwise
estimates
$$
\big|\de_\la^i\de_\xi^j\nu(\la,\xi)\big|\le
C_{i,j}|\la|^{-i}\xi^{-j}\tag6.8
$$
for $i=0,1$ and $j$ arbitrary. This implies that they can be  
appropriately
extended to the upper half-plane so that (6.4) holds with $\sigma=1$ and
$\rho$ arbitrary, so that the proposition is a consequence of Theorem 6.2.

Using the identity $\de_\xi a= \de_\xi q^\eps_\del=\frac1{2a}$, we
see by induction that
$$
\de_\xi^j\nu_0=\sum_{k=0}^j c_{jk}  
(q^\eps_\del)^{\half-k}a^{-\half+k-2j}=\nu_0\sum_{k=0}^j c_{jk}  
(q^\eps_\del)^{-k}a^{k-2j}\ .
$$

To prove (6.8) we shall use the following  elementary relation
$$\sqrt{\al+\beta}-\sqrt{\al}\,\simeq \, 
\left\{
\aligned
\frac \beta{\sqrt{\al}},\quad &\text{if}\ \al\gtrsim \beta,\\
\sqrt{\beta},\quad &\text{if}\ \al\lesssim \beta ,
\endaligned 
\right.\tag 6.9
$$ 
valid for every $\al,\beta\ge 0.$

We first consider the case $i=0$ in (6.8). Clearly, $\nu_0$ is bounded.
Moreover, 
$$a\,\simeq \, 
\left\{
\aligned
|\la|+\frac n2,\quad &\text{if}\ |\la|+\frac n2\gtrsim \sqrt{\xi},\\
\sqrt{\xi},\quad &\text{if}\ |\la|+\frac n2\lesssim \sqrt{\xi},
\endaligned 
\right.\tag 6.10
$$
and, by (6.9), since 
$$
q^\eps_\del\ge\sqrt{\xi+\la^2+\frac{n^2}4}-|\la|-\frac n2=
\sqrt{\Big(|\la|+\frac n2\Big)^2+(\xi-n|\la|)}-\sqrt{\Big(|\la|+\frac
n2\Big)^2}\  ,
$$
we have 
$$q^\eps_\del\,\gtrsim \, 
\left\{
\aligned
\frac{\xi-n|\la|}{|\la|+\frac n2},\quad &\text{if}\ |\la|+\frac n2\gtrsim
\sqrt{\xi-n|\la|},\\
\sqrt{\xi-n|\la|},\quad &\text{if}\ |\la|+\frac n2\lesssim
\sqrt{\xi-n|\la|}.
\endaligned 
\right.
$$
Notice that   $\xi-n|\la|\simeq \xi,$ since we assume 
$\xi\ge(n+1)|\la|,$ and thus 
$$q^\eps_\del\,\gtrsim \, 
\left\{
\aligned
\frac{\xi}{|\la|+\frac n2},\quad &\text{if}\ |\la|+\frac n2\gtrsim
\sqrt{\xi},\\
\sqrt{\xi},\quad &\text{if}\ |\la|+\frac n2\lesssim
\sqrt{\xi}.
\endaligned 
\right. \tag 6.11
$$
For simplicity of notation, let us assume that $\la>0.$ Then, by (6.10),
(6.11), 
$$\aligned (q^\eps_\del)^{-k}a^{k-2j}\,&\lesssim \, 
\left\{
\aligned
(\la+\frac n2)^{2k-2j}\xi^{-k},\quad &\text{if}\ \la+\frac n2\gtrsim
\sqrt{\xi},\\
(\sqrt{\xi})^{-2j},\quad &\text{if}\ \la+\frac n2\lesssim
\sqrt{\xi}
\endaligned 
\right.\\
&\lesssim \xi^{-j}. 
\endaligned
$$ 
This shows that (6.8) holds for $i=0.$

Consider now the case $i=1.$ To control $\de_\la\nu_0,$ we write
$\nu_0=\psi^\half,$ where 
$$\psi=q^\eps_\del a\inv=1+(\eps\tfrac n2+\delta\la)a\inv.$$
Then
$$\align
\de_\la\psi&=\del a\inv-(\eps\tfrac n2+\delta\la)\la a^{-3}\\
&=(\del(\xi+\tfrac {n^2}4)-\eps\tfrac n2\la)a^{-3},
\endalign
$$
hence
$$\de_\la\nu_0=\tfrac 12(\del\xi+\del\tfrac{n^2}4-\eps\tfrac
n2\la)(q^\eps_\del)^{-\half}a^{\half-3}.
$$
By induction, one then finds that 
$$\de_\xi^j\de_\la\nu_0=\nu_0\sum_{k=0}^{j-1}c_{jk}(q^\eps_\del)^{-k-1}
a^{k-2j}
+\nu_0\sum_{k=0}^{j}d_{jk}(\del\xi+\del\tfrac{n^2}4
-\eps\tfrac n2\la)(q^\eps_\del)^{-k-1} a^{k-2(j+1)}.
$$

By (6.10), (6.11) (assuming again that $\la>0$) we see that 
$$(q^\eps_\del)^{-k-1}a^{k-2j}\, \lesssim \, 
\left\{
\aligned
\frac 1{\la+\tfrac n2}(\la+\tfrac n2)^{2k-2(j-1)}\xi^{-k-1},\quad
&\text{if}\ \la+\frac n2\gtrsim \sqrt{\xi},\\
(\sqrt{\xi})^{-2j-1}\qquad,\quad &\text{if}\ \la+\frac n2\lesssim
\sqrt{\xi}.
\endaligned 
\right.
$$
Since $k\le j-1,$ these terms are of order $O((\la+\tfrac
n2)\inv\xi^{-j}).$

Noticing that $|\del\xi+\del\tfrac{n^2}4-\eps\tfrac
n2\la|\lesssim 1+\xi,$ when $\xi>(n+1)\la,$ one finds in a similar way
that the terms in the second sum are of the order 
$$(1+\xi)(\la+\tfrac n2)\inv \xi^{-(j+1)}\lesssim \la\inv \xi^{-j},
$$
so that (6.8) also holds for $i=1.$

Next, $$
\nu_+=\nu_0\sqrt{\frac a{a+\frac n2}}\ ,
$$
where the square root only depends on $\xi+\la^2$, and is a
Mihlin-H\"ormander multiplier in this variable. So Theorem 6.4 applies  
to this factor.

There remain the multipliers of the form 
$$\nu_-=(q^-_\del)^\half \big(a-\textstyle{\frac  
n2}\big)^{-\half}.$$
We begin with $i=0$ in (6.8), assuming again for
simplicity that $\la>0.$  By induction, we here see that 
$$
\de_\xi^j\nu_-=\nu_-\sum_{k+l\le j} c_{jkl}  
(q^-_\del)^{-k}(a-\tfrac n2)^{-l}a^{k+l-2j}\ .
$$
Moreover, by (6.9), 
$$\aligned
a-\frac n2 &=\sqrt{\tfrac {n^2}4+(\xi+\la^2)}-\sqrt{\tfrac
{n^2}4}\\
&\simeq \left\{
\aligned\frac{\xi+\la^2}{n/2},\quad &\text{if}\ \tfrac {n^2}4\gtrsim
\xi+\la^2,\\
\sqrt{\xi+\la^2},\quad &\text{if}\ \tfrac {n^2}4\lesssim
\xi+\la^2.
\endaligned
\right.
\endaligned
$$
Assume first that $\xi\lesssim 1$. Since $\xi>(n+1)\la,$ then 
$$a-\frac
n2\simeq \xi,\quad a\simeq 1\quad \text {and}\quad q^-_\del\gtrsim \xi
\ .
$$
Hence 
$$(q^-_\del)^{-k}(a-\tfrac n2)^{-l}a^{k+l-2j}\lesssim \xi^{-k-l}
\lesssim \xi^{-j}, 
$$
so that (6.8) holds.

 Let next $\xi\gtrsim 1$. Then $a-\frac n2\simeq \sqrt{\xi+\la^2},$
 hence  
$$
a-\frac n2 
\simeq \left\{
\aligned
\la,\quad &\text{if}\ \la\gtrsim \sqrt \xi,\\
\sqrt \xi,\quad &\text{if}\ \la\lesssim\sqrt \xi.
\endaligned
\right.
$$
In combination with (6.10), (6.11), this gives 
$$\aligned 
(q^-_\del)^{-k}(a-\tfrac n2)^{-l}a^{k+l-2j} &\lesssim  \left\{
\aligned
\Big(\frac \xi\la\Big)^{-k} \la^{-l} \la^{k+l-2j},\quad &\text{if}\
\la\gtrsim
\sqrt \xi,\\
{\sqrt \xi}^{\,-k} {\sqrt \xi}^{\,-l}{\sqrt \xi}^{\,k+l-2j},\quad
&\text{if}\
\la\lesssim\sqrt
\xi
\endaligned
\right.\\
&\lesssim \xi^{-j},
\endaligned
$$ so that (6.8) holds for $i=0.$

 Let next $i=1.$ Arguing similarly as for $\nu_0,$ we here put 
$$\psi =q^-_\del(a-\frac n2)\inv=1+\delta\la (a-\frac n2)\inv,$$
so that 
$$\de_\la\nu_-=\frac \del 2\Big(\xi a\inv (a-\frac n2)\inv-\tfrac n2
a\inv\Big) (q^-_\del)^{-\half} (a-\frac n2)^{-\half}
$$
consists of terms 
$$\mu_1=\xi (q^-_\del)^{-\half}(a-\frac n2)^{-\frac 32} a\inv, \ 
\mu_2= (q^-_\del)^{-\half}(a-\frac n2)^{-\frac 12} a\inv.
$$
Then 
$$\aligned 
\de_\xi^j\mu_1 
&=\nu_-\sum_{k+l\le
j-1}c_{jkl}\, (q^-_\del)^{-k-1}(a-\frac n2)^{-l-1} a^{-1+k+l-2(j-1)}\\
&\quad + \nu_-\sum_{k+l\le j}d_{jkl}\, \xi \, (q^-_\del)^{-k-1}(a-\frac
n2)^{-l-1} a^{-1+k+l-2j}
\endaligned
$$
and
$$ 
\de_\xi^j\mu_2 =\nu_-\sum_{k+l\le j}b_{jkl}\, (q^-_\del)^{-k-1}(a-\frac
n2)^{-l} a^{-1+k+l-2j}.
$$

If $\xi\lesssim 1,$ in view of the previous discussion one easily finds
that each term arizing in these sums is of order $O(\xi^{-1-j}),$ so that 
$|\de_\xi^j \mu_{1/2}|\lesssim \xi^{-1-j}\lesssim \la\inv\xi^{-j}.$

Similarly, if $\xi\gtrsim 1,$ then, e.g.
$$ \aligned
 \xi \, (q^-_\del)^{-k-1}&(a-\frac n2)^{-l-1} a^{-1+k+l-2j}\\
&\lesssim  \left\{
\aligned
\xi \, \Big(\frac \xi\la\Big)^{-k-1} \la^{-l-1} \la^{-1+k+l-2j},\quad
&\text{if}\
\la\gtrsim
\sqrt \xi,\\
\xi\, {\sqrt \xi}^{\,-k-1} {\sqrt \xi}^{\,-l-1}{\sqrt
\xi}^{\,-1+k+l-2j},\quad &\text{if}\
\la\lesssim\sqrt
\xi
\endaligned
\right.\\
&\lesssim \la\inv\xi^{-j}
\endaligned 
$$
if $k+l\le j,$ and the other terms can be estimated in a similar way. 

We thus see that (6.8) also holds for $i=1.$ \endproof

We can now prove our main result.

\proclaim{Theorem 6.8} If $m$ is a Mihlin-H\"ormander multiplier of  
order
$\tau>n+\half$, then $m(\Delta_1)$ is a bounded operator on
$L^p\Lambda^1(H_n)$ for $1<p<\infty$.
\endproclaim
  \proof We show that, if $m$ is as stated, then each individual
operator appearing in (6.1) is $L^p$-bounded.
We begin with the orthogonal projections and the intertwining  
operators. As
in the proof of Lemma 4.1, we write the components of  
$R=d\Delta_0^{-\half}$ as $$
(B_jL^{-\half})(L^\half\Delta_0^{-\half})\ ,\qquad (\bar
B_jL^{-\half})(L^\half\Delta_0^{-\half})\ ,\qquad T\Delta_0^{-\half}\
. $$

The Riesz transforms $ B_jL^{-\half}$, $\bar B_jL^{-\half}$ are
bounded on $L^p$, being homogeneous singular integral operators with  
smooth kernels away from
the origin. The operators
$L^\half\Delta_0^{-\half}$ and $T\Delta_0^{-\half}$ are also bounded
on $L^p$ by Theorem 6.2. In fact their spectral multipliers satisfy the  
stronger pointwise
condition (6.8) for  every $i,j$. By
duality, $R^*$ is also $L^p$-bounded.
By Corollary 4.8, $L^p$-boundedness of $P_2^+$ reduces to
$L^p$-boundedness of $\Ri$, i.e. of each operator  $$
B_j\square^{-\half}(I-\bar C)=(B_jL^{-\half})L^\half  
\square^{-\half}(I-\bar C)\ .
$$

Being homogeneous of degree zero, the spectral multiplier of $L^\half
\square^{-\half}(I-\bar C)$ satisfies (6.8), and we can apply again
Theorem 6.2. The argument is completely analogous for $P_2^-$.

It remains to discuss the last two terms. Since $\Gamma$ and
$\Gamma^*$ only contain $\Ri,\bar\Ri$ and their adjoints, we can pass  
directly to $S_\pm$ and $S^*_\pm,$ and these operators are $L^p$-bounded
by Proposition 6.7.

We finally consider  the terms containing the multiplier $m$. By Theorem
6.4,
$m(\Delta_0)$ is bounded on $L^p$, and the same is true for  
$m(\Delta_0\pm iT)$ as long
as $n\ge2$.
On $H_1$, the restriction of $m(\Delta_0\mp iT)$ to $V_2^\pm$ equals
$m(-T^2)$, as we have already observed in Section~4. Hence this case is  
even simpler,
the conclusion following by transference from $\R$ to $H_1$ (or by  
Theorem 6.2).
Finally, once we have observed that $U$ is $L^p$-bounded, it remains
to consider $m\Big(\Delta_0+\frac n2\pm
\sqrt{\Delta_0+\frac{n^2}4}\Big)$. The $L^p$-boundedness of these
operators follows from the fact that also
$$
m_\pm(s)= m\Big(s+\frac n2\pm\sqrt{s+\frac{n^2}4}\Big)
$$
satisfy a Mihlin-H\"ormander condition of order $\tau$, as a
consequence of the following last two lemmas.
\endproof

\proclaim{Lemma 6.9} Let $m$ be a Mihlin-H\"ormander multiplier of
order $\tau$ on $\R^*_+$, and let $\ph: \R^*_+\rightarrow \R^*_+$ be a  
smooth increasing function
with the following properties   \roster \item"(i)"  there exist  
exponents $\gamma$ and $\gamma'$ such that, if $k$
is the smallest integer greater than or equal to $\tau$ and $j\le k$,  
then $|\ph^{(j)}(s)|\le M s^{\gamma-j}$ for $s$ close to $0$ and  
$|\ph^{(j)}(s)|\le M
s^{\gamma'-j}$ for $s$ close to $+\infty$;
\item"(ii)" there is $\delta>0$ such that, for $j=0,1$, $\ph^{(j)}(s)\ge
\delta s^{\gamma-j}$ for $s$ close to $0$, and $\ph^{(j)}(s)\ge \delta
s^{\gamma'-j}$ for $s$ close to $+\infty$. \endroster Then $m\circ\ph$  
is also a
Mihlin-H\"ormander multiplier of order $\tau$. \endproclaim

\proof Let $I$, $J$ be compact intervals contained in $\R^*_+$, let
$\psi:I\rightarrow J$ be a $C^k$-map with never vanishing derivative,
and let $f\in L^2_\tau$ be supported on $J$. If $\tau\le k$ is an
integer, then $f\circ\psi\in L^2_\tau$, and  $$
\|f\circ\psi\|_{L^2_\tau}\le C\|f\|_{L^2_\tau}\ ,\tag6.12
$$
with $C$ depending on the $C^k$-norm of $\psi$ and on the infimum of  
$|\psi'|$.

By complex interpolation, the same is true for every $\tau\le k$.
Take now $m\in L^2_{\tau,\text{\rm sloc}}(\R^*_+)$, and consider, for  
$r>0$, $$
m_r(s)=m\circ\ph(rs)\eta_0(s)\ ,
$$
with $\eta_0$ a bump function supported in $I=[1,2]$. Define, for  
$r\le1$,
$$
\psi_r(s)=r^{-\gamma}\ph(rs)\ .
$$

By (i) and (ii), $\psi_r(I)\subseteq[\delta,M2^\gamma]=J$. The  
$C^k$-norms of the
$\psi_r$ are uniformly bounded by (i), and the derivatives $\psi'_r$
are uniformly bounded from below by (ii). By (6.12), $$  
\|m_r\|_{L^2_\tau}\le
C\|m(r^\gamma\cdot)\eta_0\circ\psi_r\inv\|_{L^2_\tau}\ . $$

Consider the set $\Omega$ consisting of the bump functions $$
\eta_r(u)= \eta_0\circ\psi_r\inv(u)=\eta_0\big(r\inv\ph\inv(r^\gamma  
u)\big)\ ,
$$
supported on $J$. It follows from (i) and (ii) that
$$
\big|(\ph\inv)^{(j)}(u)\big|\le M'u^{\frac1\gamma-j}\ ,
$$
for $j\le k$, and
$$
(\ph\inv)^{(j)}(u)\ge\delta' u^{\frac1\gamma-j}\ ,
$$
for $j=0,1$. These inequalities imply that the $\eta_r$ have uniformly
bounded $C^k$-norms. Applying now Lemma 6.1, we obtain that $$
\sup_{r\le1}\|m_r\|_{L^2_\tau}\le C\|m\|_{L^2_{\tau,\text{\rm sloc}}}\ .
$$

The same argument works for $r>1$, replacing $\gamma$ with $\gamma'$.
\endproof

\proclaim{Lemma 6.10} If $m(s)$ is a Mihlin-H\"ormander multiplier on  
$\R^*_+$ of
order $\tau>\half$, the same is true for $\tilde m(s)=m(s+a)$, for  
every $a>0$.
\endproclaim

\proof By scale-invariance, we can assume that $a=1$. Take a bump  
function $\eta_0$ with sufficiently small support, and
consider first $r$ large. By translation-invariance, the  
$L^2_\tau$-norm of $\tilde m(rs)\eta_0(s)$
equals the
$L^2_\tau$-norm of $m(rs) \eta_0(s-r\inv)$. The functions
$\eta_r(s)=\eta_0(s-r\inv)$ are supported on the same compact subset of  
$\R^*_+$, so that we can apply
Lemma 6.1 to conclude that $$ \sup_{r\ge1}\|\tilde  
m(r\cdot)\eta_0\|_{L^2_\tau}\le
C\|m\|_{L^2_{\tau,\text{\rm sloc}}}\ . $$

If we now restrict our attention to $r$ small, we can replace $m$ by
$m\chi$, where $\chi$ is smooth and supported on some interval
$[1-\delta,1+M]$. Hence we can assume that $m\in L^2_\tau$, so that
$\tilde m$ is the restriction to $\R^*_+$ of a function in $L^2_\tau$,
supported on $[-\delta,M]$. We prove that, for $r$ small, $$
\|\tilde m(r\cdot)\eta\|_{L^2_\tau}\le C\|\tilde m\|_{L^2_\tau}\  
,\tag6.13
$$
with $C$ independent of $r$.

If $\tau=k$ is an integer, it follows from Leibniz's rule that the
left-hand side is controlled by the $L^2$-norms of $r^j\tilde
m^{(j)}(rs)$ over the support of $\eta$. For $j=0$, such norms are
uniformly bounded by the boundedness of $\tilde m$, and for $j\ge1$ by
change of variable in the $L^2$-integral.
For general $\tau$, (6.13) follows by complex interpolation. \endproof


\vskip.6cm



\Refs

\widestnumber\key{MRS2}

\ref
\key AD
\by F. Astengo, B. Di Blasio
\paper The Gelfand transform of homo\-geneous distribu\-tions on
Heisen\-berg (-type) groups \jour preprint
\vol \yr \pages \endref

\ref
\key DT
\by P. de Bartolomeis, A. Tomassini
\paper On formality of some symplectic manifolds
\jour Int. Math. Res. Notices
\vol 24
\yr 2001
\pages 1287-1314
\endref

\ref
\key FS
\by G. Folland, E. M. Stein
\paper Estimates for the $\bar\de_b$-complex and analysis on the  
Heisenberg group
\jour Comm. Pure Appl. Math.
\vol 27
\yr 1974
\pages 429-522
\endref

\ref
\key G
\by D. Geller
\paper Local solvability and homogeneous distributions on the  
Heisenberg group
\jour Comm. PDE
\vol 5
\yr 1980
\pages 475-560
\endref

\ref
\key L
\by J. Lott
\paper Heat kernels on covering spaces and topological invariants
\jour J. Diff. Geom.
\vol 35
\yr 1992
\pages 471-510
\endref

\ref
\key MRS1
\by D. M\"uller, F. Ricci, E. M. Stein
\paper Marcinkiewicz multipliers and multi-parameter structure on  
Heisenberg(-type) groups, I
\jour Inv. Math.
\vol 119
\yr 1995
\pages 199-233
\endref

\ref
\key MRS2
\by D. M\"uller, F. Ricci, E. M. Stein
\paper Marcinkiewicz multipliers and multi-parameter structure on  
Heisenberg(-type) groups, II
\jour Math. Z.
\vol 221
\yr 1996
\pages 267-291
\endref

\ref
\key MS
\by D. M\"uller, E. M. Stein
\paper On spectral multipliers for Heisenberg and related groups
\jour J. Math. Pures Appl.
\vol 73
\yr 1994
\pages 413-440
\endref

\ref
\key NRS
\by A. Nagel, F. Ricci, E. M. Stein
\paper Harmonic analysis and fundamental solutions on nilpotent Lie  
groups
\jour In: Analysis and Partial Differential Equations (C. Sadosky ed.).  
Lecture Notes in
Pure Appl. Math. M. Dekker NY \vol 122
\yr 1990
\pages 249-275
\endref

\ref
\key R1
\by M. Rumin
\paper Formes diff\'erentielles sur les vari\'et\'es de contact
\jour J. Diff. Geom.
\vol 39
\yr 1994
\pages 281-330
\endref

\ref
\key R2
\by M. Rumin
\paper Sub-Riemannian limit of the differential form spectrum of  
contact manifolds
\jour GAFA
\vol 10
\yr 2000
\pages 407-452
\endref

\ref \key S
\by E. M. Stein
\book Harmonic Analysis. Real Variable Methods, Orthogonality, and  
Oscillatory Integrals
\publ Princeton Univ. Press
\yr 1993
\endref

\ref \key W
\by A. Weil
\book Introduction \`a l'\'etude des vari\'et\'es k\"ahl\'eriennes
\publ Hermann
\yr 1958
\endref

\endRefs
\enddocument